\documentclass[12pt]{amsart}
\usepackage{amsfonts}
\usepackage{amsmath}
\usepackage{amssymb}
\usepackage{url,epsfig,eqngroup}
\usepackage{color}
\usepackage{graphicx}
\usepackage{subfigure}
\usepackage{tikz}

\def\real{\mathbb{R}}
\def\cmplx{\mathbb{C}}
\def\ganz{\mathbb{Z}}

\def\eps{\varepsilon}

\newcommand\cF{\mathcal{F}}

\newcommand\cN{\mathcal{N}}
\newcommand\cO{\mathcal{O}}

\renewcommand{\Re}{\operatorname{Re}}
\renewcommand{\Im}{\operatorname{Im}}

\newcommand{\smalf}{\par\smallskip\noindent}
\newcommand{\medlf}{\par\medskip\noindent}
\newcommand{\biglf}{\par\bigskip\noindent}

\newcommand{\be}{\begin{eqnarray}}
\newcommand{\ben}{\begin{eqnarray*}}
\newcommand{\en}{\end{eqnarray}}
\newcommand{\enn}{\end{eqnarray*}}
\newcommand{\Z}{{\mathbb Z}}
\newcommand{\N}{{\mathbb N}}
\newcommand{\C}{{\mathbb C}}
\newcommand{\R}{{\mathbb R}}
\newcommand{\G}{{\Gamma}}

\newtheorem{theorem}{Theorem}[section]
\newtheorem{definition}[theorem]{Definition} 
\newtheorem{lemma}[theorem]{Lemma}
\newtheorem{corollary}[theorem]{Corollary}
\newtheorem{proposition}[theorem]{Proposition}
\newtheorem{remark}[theorem]{Remark}

\newtheorem{assumption}[theorem]{Assumption}

\definecolor{rot}{rgb}{0,0,0}

\definecolor{rot1}{rgb}{1,0,0}

\begin{document}

\title{Direct and inverse time-harmonic scattering by Dirichlet periodic curves with local 
perturbations}

\author{Guanghui Hu}
\address{Guanghui Hu: School of Mathematical Sciences and LPMC\\
Nankai University \\
Tianjin 300071, China}
\email{ghhu@nankai.edu.cn}
\author{Andreas Kirsch}
\address{Andreas Kirsch: Department of Mathematics \\
Karlsruhe Institute of Technology (KIT) \\
76131 Karlsruhe, Germany}
\email{andreas.kirsch@kit.edu}

\date{\today}

\begin{abstract} This is a continuation of the authors' previous work (A. Kirsch, Math. Meth. Appl. Sci., 45 (2022): 5737-5773.) on well-posedness of time-harmonic scattering by locally perturbed periodic curves of Dirichlet kind. The scattering interface is supposed to be given by a non-self-intersecting Lipschitz curve. We study properties of the Green's function and prove new well-posedness results for scattering of plane waves at a propagative wave number. In such a case there exist guided waves to the unperturbed problem, which are also known as  Bounded States in the Continuity (BICs) in physics.  
In this paper uniqueness of the forward scattering follows from an orthogonal constraint condition enforcing on the total field to the unperturbed scattering problem. This constraint condition, which is also valid under the Neumann boundary condition, is derived from the singular perturbation arguments and also from the approach of approximating a plane wave by point source waves. For the inverse problem of determining the defect, we prove several uniqueness results using a finite or infinite number of point source and plane waves, depending on whether a priori information on the size and height of the defect is available. 

\vspace{.2in} 

Keywords: Helmholtz equation, non-self-intersecting periodic curve, local perturbation, Dirichlet boundary condition, plane wave, uniqueness, inverse problem.
\end{abstract}

\maketitle

\section{Introduction}

This paper is concerned with the TE polarization of time-harmonic electromagnetic scattering from perfectly conducting gratings with a localized defect. The first part deals with 
well-posedness of the mathematical model for plane wave incidences and properties of the Green's function. In the second part, we study uniqueness to the inverse  problems of determining the local perturbation from near/far-field data excited by plane and point source waves. Throughout the paper, the cross-section of the scattering surface is supposed to be a  non-self-intersecting periodic curve with a local perturbation. In the TE polarization case, the grating diffraction problem can be modeled by the Dirichlet boundary value problem of the two-dimensional Helmholtz equation in the unbounded domain above the interface complemented with a proper radiation radiation at infinity. We refer to \cite{Bao2022,P1980} for a comprehensive introduction of electromagnetic scattering theory for diffraction gratings. 
 
\smalf

At the absence of the defect, the wave field for a plane wave incidence is well-known to be quasiperiodic, due to the periodicity of the scattering interface and the quasi-periodicity of the incoming plane wave. The Rayleigh expansion radiation condition (which was originally proposed by Lord Rayleigh in 1907 \cite{R07}) has been widely used in the literature concerning the mathematical analysis and numerical approximation of wave scattering in periodic structures. With the Fredholm theory, it is also well known that the forward scattering model is well-posed for all incident frequencies excluding a discrete set with the only accumulating point at infinity. However, the Rayleigh expansion radiation condition does not always lead to uniqueness (although existence can always be justified via variational argument), because of the existence of guided/Floquet wave modes to the homogeneous problem, which exponentially decay in the direction orthogonal to the periodicity direction \cite{BBS94, G00, HR15, KN02}.  
If the interface is given by the graph of some periodic function (a weaker condition was proposed in \cite{CE10, SM05}), uniqueness and existence can be proved (which implies the absence of guided waves) for the Dirichlet boundary value problem of the Helmholtz equation at an arbitrary frequency; see \cite{EY02, K93}.
Since an incoming point source is not quasi-periodic, the Rayleigh expansion condition is not valid any more. Instead, the Upward Propagation Radiation Condition \cite{CZ98} or the
Angular Spectrum Representation Condition \cite{SM05, CE10} can be used for proving well posedness within the framework of rough surface scattering problems, provided the domain with a geometrical condition admits no guided waves.  Since a locally perturbed periodic curve  can be treated as a special rough curve, it was proved in \cite{HWR} that the Green's function to the perturbed scattering problem satisfies a Half-plane Sommerfeld radiation condition and the scattered field generated by a plane wave and caused by the defect fulfills the same radiation condition, as long as guided modes can be excluded.   

\smalf
The mathematical analysis is more involved for locally perturbed scattering problems when guided waves exist in periodic structures. An open wave-guide radiation condition (which is equivalent to the closed wave-guide radiation condition \cite{FJ16} based on dispersion curves) was proposed in \cite{KL-MMAS} for acoustic scattering by inhomogeneous periodic layers in a half-plane. Such a radiation condition was derived from the Limiting Absorption Principle and the Floquet-Bloch transform and was later extended to investigate well-posedness of wave scattering by layered periodic media in $\R^2$ and by periodic tubes in $\R^3$ (see \cite{TF, K19-1, K19-2, KL18, K22, KL-MMAS}). This open wave guide radiation 
condition consists of a radiating part and a 
propagating (guided) part. It was recently shown in \cite{K22} that the radiating part satisfies a Sommerfeld-type radiation condition and, due to the existence of cut-off values, the radiating part  decays as $|x_1|^{-1/2}$ in the periodicity direction. 
In the authors' previous work \cite{HK22}, the open wave guide radiation condition has been adopted to prove well-posedness of Dirichlet and Neumann boundary value problems of the Helmholtz equation in a locally perturbed periodic structure. 
By constructing a Dirichlet-to-Neumann operator on the boundary of a truncated domain,
 uniqueness and existence of time-harmonic scattering by incoming point 
source waves, plane waves and surface waves are 
established.  This has generalized the results of \cite{HWR} to scattering interfaces given by non-self-intersecting curves, for which the forward solutions may contain guided wave modes.

\smalf

For plane wave incidence, the well-posedness results of \cite{HK22}  are based on the uniqueness assumption on the forward scattering model in periodic structures.  This is equivalent to the statement that the quasi-periodicity of the incoming plane wave (that is $k\sin\theta$, where $\theta$ is the incident angle) is not a propagative wave-numbmer (see Definition \ref{d-exceptional} (ii)).  If otherwise, solutions to the unperturbed scattering problem are not unique and neither for the perturbed problem; see Section \ref{sec3-1} for detailed discussions. In the first part of this paper, we shall propose an additional constraint condition on solutions of the unperturbed problem to ensure uniqueness. For this purpose we adopt two different approaches to the
Dirichlet boundary value problem: the Limiting Absorption Principle by replacing $k$ 
by $k+i\epsilon$ (Subsection~\ref{LAP}) and the approximation by point source waves 
(Subsection~\ref{APS}). It will be shown in Theorems \ref{TH-LAP} (i) and \ref{TH-POI} that both methods yield the same constraint condition. The limiting absorption arguments for approximating wave numbers has complemented the work of \cite{K22p}, where the LAP for approximating refractive indices, the continuity with respect to incident angles together with the method of approximating plane waves by point source waves were
justified for scattering by layered periodic media. 
In the first part we also justify some properties of the Green's function to perturbed and unperturbed scattering problems; see Sections \ref{sec3} and \ref{sec3-1}. In particular, the mixed reciprocity relation between point source and plane wave incidences will be verified in Theorem \ref{TH-PL}. 

\smalf

We remark that radiation conditions and numerical approximations with exact boundary conditions (DtN maps) were also considered in \cite{FJ09, JLF2006}
for wave propagating in a closed 
periodic wave-guide  and in a photonic crystal  containing a local 
perturbation. We refer to \cite{LZ17, ZR} for numerical methods based on  LAP and the Floquet-Bloch transform and to \cite{XHLR} using the boundary integral equation method in combination with perfectly mathched absorbing layers.  

\smalf
The second part of this paper concerns  inverse scattering problems of recovering the localized defect by assuming {\rm a priori} knowledge on the unperturbed periodic structure.  Using infinitely many point source or plane waves at a fixed energy, we prove that   
the position and shape of the local defect can be uniquely determined by the corresponding near-field data measured on a line segment above the interface; see Subsection \ref{sub5-1}. 
As will be seen in the proof of Theorems \ref{uni-plane-inf}, the complexity of the solution structure gives rise to essential difficulties in justifying linear independence of the wave fields for different angles. If some a priori information on the defect is available, one can prove uniqueness with a finite number of incoming waves by adopting Colton and Slemann's idea of determining a bounded sound-soft obstacle \cite{CS83} (see also \cite{HF97} for the corresponding results in periodic structures with a fixed direction). A counterexample will be constructed to show that one incident plane wave is impossible to imply uniqueness in general. 
In Subsection \ref{sub5-4}, we discuss uniqueness results using far-field patterns of incoming point source waves over finite or infinite number of observation directions. 

In all of the paper, we choose the square root function to be 
holomorphic in the cutted plane $\cmplx\setminus(i\real_{\leq 0})$. In particular, 
$\sqrt{t}=i\sqrt{|t|}$ for $t\in\real_{<0}$.  The functions $\phi$ are called guided (or 
propagating or Floquet) modes.
The Fourier transform is defined as
$$ (\cF\phi)(\omega)\ :=\ \frac{1}{\sqrt{2\pi}}\int\limits_{-\infty}^\infty\phi(s)\,
e^{-is\omega}\,ds\,,\quad\omega\in\real\,, $$
which can be considered as an unitary operator from $L^2(\real)$ onto itself. For a domain $\Omega\subset \R^2$,
the weighted Sobolev space $H_\rho^1(\Omega)$ is 
defined by
\ben
H_\rho^1(\Omega)\ :=\ \bigl\{u:(|1+|x_1|^2)^{\rho/2}u\in H^1(\Omega)\bigr\}\,,
\quad \rho\in\R\,.
\enn

\section{Radiation 
conditions and well-posedness results}
\label{s-intro}

In this section we describe the mathematical model for the TE polarization of time-harmonic electromagnetic scattering from a perfectly conducting periodic surface with local perturbations. We first recall some notations, define the open wave-guide radiation and then present some well-posedness results from the authors' previous work \cite{HK22}. 

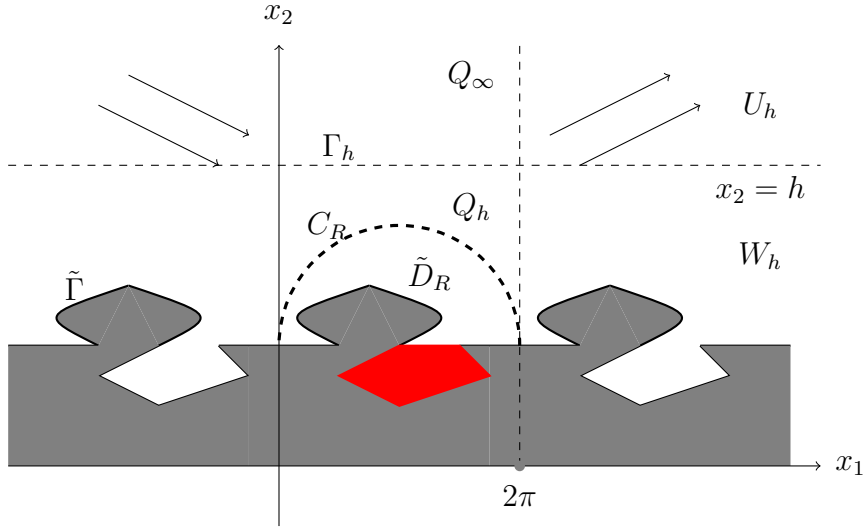
\begin{figure}[htb]
  \centering
  \begin{tikzpicture}[scale=.8, thick]
\draw (0.5,0) -- (2,0);
\filldraw[fill=gray] (2,0) .. controls (1, 0.5) .. (2.5,1);
\filldraw[fill=gray] (2.5,1) .. controls (4, 0.5) .. (3,0);
\draw (3, 0) -- (2, -0.5) -- (3, -1) -- (4.5, -0.5) -- (4, 0);

\draw (4,0) -- (6,0);
\filldraw[fill=gray] (6,0) .. controls (5, 0.5) .. (6.5,1);
\filldraw[fill=gray] (6.5,1) .. controls (8, 0.5) .. (7,0);
\draw (7, 0) -- (6, -0.5) -- (7, -1) -- (8.5, -0.5) -- (8, 0);


\draw (8,0) -- (10,0);
\filldraw[fill=gray] (10,0) .. controls (9, 0.5) .. (10.5,1);
\filldraw[fill=gray] (10.5,1) .. controls (12, 0.5) .. (11,0);
\draw (11, 0) -- (10, -0.5) -- (11, -1) -- (12.5, -0.5) -- (12, 0);
\draw (12,0) -- (13.5,0);
\fill[gray] (0.5,-2) -- (0.5,0) -- (2,0) -- (2.5,1) -- (3,0) -- (2, -0.5) -- (3, -1) -- (4.5, -0.5) -- (4, 0)--(4.5, 0) -- (4.5, -2);
\fill[gray] (4.5,-2) -- (4.5,0) -- (6,0) -- (6.5,1) 
-- (7,0) -- (6, -0.5) -- (7, -1) -- (8.5, -0.5) -- (8, 0)--(8.5, 0) -- (8.5, -2);
\fill[gray] (8.5,-2) -- (8.5,0) -- (10,0) -- (10.5,1) -- (11,0) -- (10, -0.5) -- (11, -1) -- (12.5, -0.5) -- (12, 0)--(13.5, 0) -- (13.5, -2);
\filldraw[red]  (7,0) -- (6, -0.5) -- (7, -1) -- (8.5, -0.5) -- (8, 0)--(7,0);


\draw[very thick, dashed] (9,0) arc (0:180:2);
\draw[thin,->] (5,-3) -- (5,5);
\draw[thin,->] (0.5,-2) -- (14,-2);

\draw[thin, dashed] (0.5,3) -- (14,3);
\draw[thin, dashed] (9,-2) -- (9,5);
\draw[thin,->] (2,4) -- (4,3);
\draw[thin,->] (2.5,4.5) -- (4.5,3.5);
\draw[thin,->] (10,3) -- (12,4);
\draw[thin,->] (9.5,3.5) -- (11.5,4.5);

\filldraw [gray] (9,-2) circle (2pt);
\draw (9,-2.5) node{$2\pi$};
\draw (14.5,-2) node{$x_1$};
\draw (5,5.5) node{$x_2$};
\draw (6,3.3) node{$\Gamma_h$};
\draw (13,2.6) node{$x_2=h$};
\draw (13,4) node{$U_h$};
\draw (13,1.5) node{$W_h$};
\draw (8.2,2.3) node{$Q_h$};
\draw (8.2,4.5) node{$Q_\infty$};
\draw (1.6,1) node{$\tilde{\Gamma}$};
\draw (5.8,2) node{$C_R$};
\draw (7.5,1.2) node{$\tilde{D}_R$};
\end{tikzpicture}
\caption{Illustration of wave scattering from a perfectly reflecting periodic curve with a local perturbation in $(0, 2\pi)$. The red area denotes the perturbed domain. The scattering interface is supposed to be a non-self-intersecting curve. }\label{f12}
\end{figure}

Let $D\subset\R^2$ be a $2\pi$-periodic domain with respect to the $x_1$-direction. 
The boundary $\Gamma:=\partial D$ is supposed to be given by a non-self-intersecting Lipschitz 
curve which is bounded in $x_2$-direction and $2\pi$-periodic with respect to $x_1$. Let 
$\tilde{D}$ be a local perturbation of $D$ in the way that $\Gamma\setminus\tilde{\Gamma}$ and 
$\tilde{\Gamma}\setminus\Gamma$ are bounded where $\tilde{\Gamma}=\partial \tilde{D}$ is 
the perturbed boundary which is also assumed to be a non-self-intersecting curve (see Figure \ref{f12}). Suppose 
that $\tilde{D}$ is filled by a homogeneous and isotropic medium and that $\tilde{\Gamma}$ 
is a perfectly reflecting curve of Dirichlet kind. Let $u^{in}$ be an incoming wave incident 
onto $\tilde{\Gamma}$. The scattered field $u^{sc}$ can be governed by the boundary value 
problem of the Helmholtz equation
\ben
\Delta u^{sc}+k^2 u^{sc}=0\quad\mbox{in}\quad \tilde{D},\qquad u^{sc}=-u^{in}\quad\mbox{on}\quad 
\tilde{\G},
\enn
complemented by some radiation condition in $\tilde{D}$ explained below. To specify this 
radiation condition we need to introduce several definitions and make some assumptions.  

For the forward scattering problem, we suppose without loss of generality (changing the 
period of the periodic structure if otherwise) that the perturbations 
$\Gamma\setminus\tilde{\Gamma}$ and $\tilde{\Gamma}\setminus\Gamma$ are contained in the disc 
$\{x\in\real^2:(x_1-\pi)^2+x_2^2<\pi^2\}$. We fix $R, h_0>\pi$ throughout this paper and 
use the following notations for $h>h_0$ (see Figure \ref{f12}).
\begin{eqnarray*}
Q_h\ & := & \{x\in D:0<x_1<2\pi,\ x_2<h\}\,,\quad Q_\infty\ :=\ \{x\in D:0<x_1<2\pi\}\,, \\
\Gamma_h & := & (0,2\pi)\times\{h\}\,,\quad
W_h\ :=\ \{x\in D:x_2<h\}\,,\quad U_h\ :=\ \{x\in D: x_2>h\}\, , \\
C_R & := & \{x\in D:(x_1-\pi)^2+x_2^2=R^2\},\quad 
\Sigma_R\ :=\ \{x\in D: (x_1-\pi)^2+x_2^2>R^2\}\,, \\
D_R & := & \{x\in D:(x_1-\pi)^2+x_2^2<R^2\}\,,\quad
\tilde{D}_R\ :=\ \{x\in \tilde{D}:(x_1-\pi)^2+x_2^2<R^2\}.
\end{eqnarray*}

We recall that a function $\phi\in L^2_{loc}(\real)$ is called $\alpha$-quasi-periodic if 
$\phi(x_1+2\pi)=e^{2\pi\alpha i}\phi(x_1)$ for all $x_1\in\real$.
Below we introduce some function 
spaces \footnote{The definitions hold also for $D$ instead of $\tilde{D}$}.
\begin{eqnarray*}
H^1_{loc}(\tilde{D}) & := & \bigl\{u|_{\tilde{D}}:u\in H^1_{loc}(\real^2)\bigr\}\,, \\
H^1_{loc,0}(\tilde{D}) & := & \bigl\{u\in H^1_{loc}(\tilde{D}):u=0\mbox{ on }
\partial\tilde{D}\bigr\}\,, \\
H^1_\ast(\Sigma_R) & := & \biggl\{u\in H^1_{loc}(\Sigma_R):
\begin{array}{l} u|_{W_h\cap\Sigma_R}\in H^1(W_h\cap\Sigma_R)\;\mbox{for all}\;h>h_0,\\
u=0\;\mbox{on}\;\partial\Sigma_R\cap\partial D \end{array}\biggr\}\,, \\
H^1_{\alpha,loc}(D) & := & \bigl\{u\in H^1_{loc}(D):u(\cdot,x_2)\mbox{ is 
$\alpha$-quasi-periodic} \bigr\}\,, \\
H^1_{\alpha,loc,0}(D) & := & \bigl\{u\in H^1_{\alpha,loc}(D): u=0\mbox{ on }\partial D
\bigr\}\,.
\end{eqnarray*}

\begin{definition} \label{d-exceptional}
(i) $\alpha\in[-1/2,1/2]$ is called a \emph{cut-off value} if there exists $\ell\in\ganz$ 
such that $|\alpha+\ell|=k$.
\newline
(ii) $\alpha\in[-1/2,1/2]$ is called a \emph{propagative wave number} if there exists a 
non-trivial $\phi\in H^1_{\alpha,loc,0}(D)$ such that 
\begin{eqngroup}\begin{equation} \label{exc:a}
\Delta\phi + k^2\phi\ =\ 0\text{ in }D\,,
\end{equation}
and $\phi$ satisfies the upward Rayleigh expansion
\begin{equation} \label{exc:b}
\phi(x)\ =\ \sum_{\ell\in\ganz}\phi_\ell\,e^{i(\ell+\alpha)x_1}\,
e^{i\sqrt{k^2-(\ell+\alpha)^2}(x_2-h_0)}\quad\mbox{for }x_2>h_0
\end{equation}\end{eqngroup}
for some $\phi_\ell\in\cmplx$ where the convergence is uniform for $x_2\geq h_0+\eps$ for 
every $\eps>0$.  
\end{definition}
\begin{remark}
In physical literatures (see e.g., \cite{Lu19, Lu18}), the guided model $\phi(\cdot, \alpha, k)$ of \eqref{exc:b} is called 
 a Bound State in the Continuity (BIC) if $|\alpha|<k$.
\end{remark}

In Definition~\ref{d-exceptional} we restrict the quasi-periodic parameter $\alpha$ to the 
interval $[-1/2, 1/2]$, because an $\alpha$-quasi-periodic function must be also 
$(\alpha+j)$-quasi-periodic for any $j\in \N$. The possible existence of guided waves leads to essential difficulties in proving well-posedness of forward scattering problems under the Rayleigh expansion condition \eqref{exc:b}, because they are solutions to the homogeneous problem when $u^{in}= 0$.
It is worthy mentioning that the set of propagative wave numbers must be empty, if the unperturbed domain $D$ fulfills the following geometrical condition (see \cite{SM05, CE10}): 
\begin{equation}\label{GCD}
(x_1,x_2)\in D\quad\Longrightarrow\quad(x_1, x_2+s)\in D\quad\mbox{for all }
s>0\,.
\end{equation}
In the special case that  
$\Gamma$ is given by the graph of some function, uniqueness was verified in \cite{K93, EY02} under different regularity assumptions enforcing on $\Gamma$. Below we discuss the existence of propagative wavenumbers. 
Throughout this paper we make the following 
assumptions.
\begin{assumption} \label{assump1}
Let $|\ell+\alpha|\not=k$ for every propagative wave number $\alpha\in[-1/2,1/2]$ and every 
$\ell\in\ganz$; that is, no cut-off value is a propagative wave number.
\end{assumption}
\smalf
Note that this assumption can be automatically fulfilled if the geometrical condition \eqref{GCD} holds true.
Under Assumption \ref{assump1} it can be shown (see, e.g. \cite{KL18} for the case of a flat 
curve $\Gamma=\Gamma_0$ and an additional index of refraction) that at most a finite number of 
propagative wave numbers exists in the interval $[-1/2,1/2]$. Furthermore, if $\alpha$ is a 
propagative wave number with mode $\phi$ then $-\alpha$ is a propagative wave number with 
mode $\overline{\phi}$. Therefore, we can number the propagative wave numbers in $[-1/2,1/2]$ 
such that they are given by $\{\hat{\alpha}_j:j\in J\}$ where $J\subset\ganz$ is finite and 
symmetric with respect to $0$ and $\hat{\alpha}_{-j}=-\hat{\alpha}_j$ for $j\in J$. 
Furthermore, it is known that (under Assumption~\ref{assump1}) every mode $\phi$ is evanescent;
that is, exponentially decaying as $x_2$ tends to infinity in $D$; that is, satisfies 
$|\phi(x)|\leq c\,e^{-\delta|x_2|}$ for $x_2\geq h_0$ and some $c,\delta>0$ which are 
independent of $x$. The corresponding space
\begin{equation} \label{X_j}
X_j\ :=\ \bigl\{ \phi\in H^1_{\hat{\alpha}_j,loc,0}(D):u\mbox{ satisfies (\ref{exc:a}) and 
(\ref{exc:b}) for }\alpha=\hat{\alpha}_j \bigr\}
\end{equation}
of modes is finite dimensional with some dimension $m_j>0$. We refer to Lemma \ref{lem3.4} for discussions on these properties in periodic Sobolev spaces. 

\smalf

On $X_j$ we define the 
sesqui-linear form $B:X_j\times X_j\to\cmplx$ by 
\begin{equation} \label{e-sesqui}
B(\phi,\psi)\ :=\ -2i\int\limits_{Q_\infty}\frac{\partial\phi}{\partial x_1}\,
\overline{\psi}\,dx\,,\quad\phi,\psi\in X_j\,.
\end{equation}
Using integration by parts and the exponential decay of $\phi\in X_j$, we obtain $B(\phi,\psi)=\overline{B(\psi,\phi)}$ for all $\phi,\psi\in X_j$. This implies that 
 $B$ is hermitian and that $B(\phi, \phi)$ is real valued for all $\phi\in X_j$.
Now we assume that $B$ is non-degenerated on every 
$X_j$ in the sense that
\begin{assumption} \label{assump2}
For every $j\in J$ and $\psi\in X_j$, $\psi\not=0$, the linear form  
$B(\cdot,\psi):X_j\to\cmplx$ is non-trivial on $X_j$; that is, there exists $\phi\in X_j$ 
with $B(\phi,\psi)\not=0$.
\end{assumption}
The hermitian sesqui-linear form $B$ defines the cones $\{\psi\in X_j:B(\psi,\psi)\gtrless 0\}$ 
of propagating waves traveling to the right and left, respectively (see also \cite{FJ16,K19-1, KL18, K22}). We construct a basis of 
$X_j$ with elements in these cones by taking the inner product $(\cdot,\cdot)_{X_j}$ and 
consider the following eigenvalue problem in $X_j$ for every fixed $j\in J$. Determine 
$\lambda_{\ell,j}\in\real$ and non-trivial $\hat{\phi}_{\ell,j}\in X_j$ with 
\begin{equation} \label{evp}
B(\hat{\phi}_{\ell,j},\psi)\ =\ -2i\int\limits_{Q_\infty}\frac{\partial\hat{\phi}_{\ell,j}}
{\partial x_1}\,\overline{\psi}\,dx\ =\ 
\lambda_{\ell,j}\,\bigl(\hat{\phi}_{\ell,j},\psi\bigr)_{X_j}\quad\mbox{for all }\psi\in X_j
\end{equation}
and $\ell=1,\ldots,m_j$. We normalize the basis such that
$\bigl(\hat{\phi}_{\ell,j},\hat{\phi}_{\ell^\prime,j}\bigr)_{X_j}=\delta_{\ell,\ell^\prime}$ 
for $\ell,\ell^\prime=1,\ldots,m_j$. Then $\lambda_{\ell,j}=B(\hat{\phi}_{\ell,j},
\hat{\phi}_{\ell,j})$ and the Assumption~\ref{assump2} is equivalent to 
$\lambda_{\ell,j}\not=0$ for all $\ell=1,\ldots,m_j$ and $j\in J$. Physically, the assumption \eqref{assump2}
is equivalent to the assumption that the group velocity of each guided mode is non-vanishing (see \cite[Remark 1.4]{K22}). 
\medlf
Now we are able to formulate the open waveguide radiation condition for our Dirichlet boundary value 
problem (see \cite{HK22}).

\begin{definition} \label{d-RC}
Let $\psi_+,\psi_-\in C^\infty(\real)$ be any functions with $\psi_\pm(x_1)=1$ for 
$\pm x_1\geq\sigma$ (for some $\sigma>\max\{R,2\pi\}+1$) and $\psi_\pm(x_1)=0$ 
for $\pm x_1\leq\sigma-1$. 
\smalf
A solution $u\in H^1_{loc}(\Sigma_R)$ of the Helmholtz equation $(\Delta+k^2)u=0$ satisfies 
the open waveguide radiation condition with respect to an inner product $(\cdot,\cdot)_{X_j}$ 
in $X_j$ if $u$ has in $\Sigma_R  $ a decomposition into $u=u_{rad}+u_{prop}$ which satisfy 
the following conditions.
\begin{itemize}
\item[(a)] The propagating part $u_{prop}$ has the form 
\begin{equation} \label{u2}
u_{prop}(x)\ =\ \sum_{j\in J}\biggl[\psi_+(x_1)\sum_{\ell:\lambda_{\ell,j}>0}
a_{\ell,j}\,\hat{\phi}_{\ell,j}(x)\ +\ \psi_-(x_1)\sum_{\ell:\lambda_{\ell,j}<0}a_{\ell,j}\,
\hat{\phi}_{\ell,j}(x)\biggr]
\end{equation}
for $x\in\Sigma_R  $ and some $a_{\ell,j}\in\cmplx$. Here, for every $j\in J$ the 
scalars $\lambda_{\ell,j}\in\real$ and $\hat{\phi}_{\ell,j}\in\hat{X}_j$ for 
$\ell=1,\ldots,m_j$ are given by the eigenvalues and corresponding eigenfunctions, 
respectively, of the self adjoint eigenvalue problem (\ref{evp}). 
\item[(b)] The radiating part $u_{rad}\in H^1_\ast(\Sigma_R)$ satisfies the generalized 
angular spectrum radiation condition
\begin{equation} \label{angulaa-spectrum-rc}
\int\limits_{-\infty}^\infty\left|\frac{\partial(\cF u_{rad})(\omega,x_2)}{\partial x_2}-
i\sqrt{k^2-\omega^2}\,(\cF u_{rad})(\omega,x_2)\right|^2d\omega\ \longrightarrow\ 0\,,\quad 
x_2\to\infty\,.
\end{equation}

\end{itemize}
\end{definition}
\smalf
The above radiation condition has been earlier studied in \cite{TF, KL18,K22, K19-1, K19-2, KL-MMAS} for layered periodic structures. 
It has been shown in \cite{KL18} for the case of 
half plane source problem with an inhomogeneous period layer that the radiation condition of 
Definition \ref{d-RC} for the inner product $(\phi,\psi)_{X_j}=2k\int\limits_{Q_\infty}n\,\phi\,
\overline{\psi}\,dx$ is a consequence of the limiting absorption principle by replacing $k$ 
with $k+i\epsilon$, $\epsilon>0$. Here the function $n$ stands for the refractive index of the inhomogeneous layer.

\begin{assumption}[Absence of bounded states]\label{assump3}
There are no bound states to the perturbed scattering problem, that is, any solution 
$u\in H^1_0(\tilde{D})$ of $\Delta u+k^2u=0$ in $\tilde{D}$ must vanish identically. 
\end{assumption}
One can remove this assumption if the domain $\tilde{D}$ fulfils the  condition \eqref{GCD}. Note that with this geometrical
condition on $\tilde{D}$, the unperturbed domain $D$ should also meet the requirement  \eqref{GCD} and thus  the existence of propagating modes $\hat{\phi}_{\ell,j}$ is excluded \cite{SM05}.
\medlf
In all of the paper we make Assumptions~\ref{assump1},  \ref{assump2} and \ref{assump3} without 
mentioning this any more. We consider two kinds of incoming waves: 
\begin{itemize}
\item[(i)] Point source wave: $u^{in}(x):=\Phi(x;y)=\frac{i}{4}H_0^{(1)}(k|x-y|)$ with the 
source position $y\in \tilde{D}_R$. 
\item[(ii)] Plane wave: $u^{in}(x)=e^{ikx\cdot\hat{\theta}}$ where 
$\hat{\theta}=(\sin\theta,-\cos\theta)$ is the incident direction with some incident angle 
$\theta\in(-\pi/2, \pi/2)$. 
\end{itemize}


Before stating uniqueness and existence results, we recall the Sommerfeld radiation condition 
used in \cite{HWR, K22}.
\begin{definition}\label{src} 
A function $v\in C^\infty(U_{h_0}\cap\Sigma_R  )$ satisfies the Sommerfeld radiation 
condition in $U_{h_0}\cap\Sigma_R$ if $v\in H^1_\rho(W_h\cap\Sigma_R)$ for all 
$h>h_0$ and all $\rho<1$ and
\be\label{SRC}
\sup\limits_{x\in C_a\cap U_h}|x|^{1/2}\,\Bigl|\frac{\partial v(x)}{\partial r}-ik v(x)\Big|\ 
\rightarrow\ 0\,,\quad a\rightarrow\infty\,,\qquad 
\sup\limits_{x\in U_h}|x|^{1/2}|v(x)|\ <\ \infty
\en
for all $h>h_0$ where $r=|x|$.
\end{definition}
It was shown in the authors' previous paper \cite{HK22} that the radiation condition for the radiating part of the open waveguide radiation condition of Def. \ref{d-RC} is equivalent to the above the Sommerfeld radiation condition.
We remark that the point source wave $\Phi(x; y)$ with $y\in D_R$ satisfies the Sommerfeld radiation condition of Def. \ref{src} with the index $\rho<0$, because $\Phi(x;y)\sim |x|^{-1/2}$ as $|x|\rightarrow\infty$ in $\R^2$. However, plane waves and quasi-periodic surface waves are not included. Such kinds of wave modes belong to  
$ H^1_\rho(W_h\cap\Sigma_R)$ for all 
$h>h_0$ with the index $\rho<-1/2$.
An integral form of the above Sommerfeld radiation condition is defined as follows. 
\begin{definition}\label{src-i} 
Let $a_j$ be a sequence in $\real$ such that $a_j\to\infty$ and suppose $\tilde{D}_{a_j}$ are 
Lipschitz domains. A solution $v\in H^1_{loc}(\Sigma_R)$ 
satisfies the Sommerfeld radiation condition in integral form if 
$$ \biggl\Vert\frac{\partial v}{\partial r}-ik v\biggr\Vert_{H^{-1/2}(C_{a_j})}\ 
\longrightarrow\ 0\,,\quad j\to\infty\,, $$
where $r=|x|$.
\end{definition}
\begin{lemma}\label{coro4.1} 
\begin{itemize}
\item[(i)] If $v$ satisfies the Sommerfeld radiation condition of 
Definition~\ref{src} with the index $\rho\geq 0$, then $v$ also fulfills the integral form of 
the radiation condition defined by Definition~\ref{src-i}.
\item[(ii)] The condition (b) for the radiating part of $u$ in Def.
\ref{d-RC} is equivalent with the Sommerfeld radiation condition of Def. \ref{src}.
\end{itemize}
\end{lemma}
Lemma \ref{coro4.1} and the following well-posedness results for incident plane and point source waves were proved 
in the authors' previous paper \cite{HK22}. 

\begin{proposition}[Well-posedness for point source waves]\label{wps}  
Let $u^{in}:=\Phi(\cdot\, ; y)$ be an incoming point source wave with $y\in\tilde{D}_R$. Then 
the locally perturbed scattering problem admits a unique solution $u$
such that $u-u^{in}\in H^1_{loc}(\tilde{D})$ and $u$ satisfies the open waveguide radiation 
condition of Definition \ref{d-RC}. Furthermore, the radiating part $u_{rad}$ of $u$ satisfies the 
Sommerfeld radiation conditions of Definitions~\ref{src} and \ref{src-i}.
\end{proposition}
If  $\partial\tilde{D}$ is given by a Lipschitz graph (that is, guided waves are excluded), 
the results of Proposition \ref{wps} were verified in \cite{HWR} within a more general framework for rough surface scattering problems. We remark that, in such a case, the scattered field $u^{sc}:=u-u^{in}$ does not satisfy the open waveguide radiation condition of Def. \ref{d-RC} and either the Sommerfeld radiation condition of Def. \ref{src}, because $u^{in}=\Phi(\cdot; y)$ does not belong to $H^1(\Sigma_R)$. In fact, $u^{sc}$ fulfills the Sommerfeld radiation condition with the index $\rho<0$.

\begin{proposition}[Well-posedness for plane waves]\label{wpw}
Let $\alpha:=k\sin\theta$ be not a propagative wave number (see Definition~\ref{d-exceptional} 
(ii)). Then the perturbed scattering problem for a plane wave incidence 
$u^{in}(x)=e^{ikx\cdot\hat{\theta}}$ admits a unique solution 
$u=u^{in}+u^{sc}\in H_{loc,0}^1(\tilde{D})$ such that the scattered part $u^{sc}$ has a 
decomposition in the form $u^{sc}=u^{sc}_{unpert}+u^{sc}_{pert}$ in the region $\Sigma_R  $
where $u^{sc}_{unpert}\in H^1_{\alpha,loc}(D)$ is the scattered field corresponding to the 
unperturbed problem that satisfies the upward Rayleigh expansion~\eqref{exc:b} 
with the quasi-periodic parameter $\alpha=k\sin\theta$. The part $u^{sc}_{pert}\in 
H_{loc}^1(\Sigma_R)$ fulfils the open waveguide radiation condition of Definition \ref{d-RC}
and the radiating part of $u^{sc}_{pert}$ satisfies the Sommerfeld radiation conditions of Definitions \ref{src} and \ref{src-i}.
\end{proposition}

We emphasize in Proposition \ref{wpw} that, since $k\sin\theta$ is assumed to be no propagative wavenumber (or equivalently, no BIC exists at the pair $(\alpha, k)\in \R^2$ with $\alpha=k\sin\theta$), the unperturbed scattered field $u^{sc}_{unpert}$ is unique. In the subsequent Sections \ref{sec3} and \ref{sec3-1}, we shall carry out further studies on forward scattering problems, including properties of the Green's functions to perturbed and unperturbed problems and well-posedness for a plane wave incidence at a propagative wavenumber (that is, when a BIC occurs).   Theorems \ref{sy},   \ref{TH-LAP} and \ref{TH-POI} will be used later for investigating inverse problems in Section \ref{sec4}.

\section{Properties of Green's Function}\label{sec3}

Let $\Phi(x; y)$ be the fundamental solution of the Helmholtz equation. 
The Green's function $G$ for $\tilde{D}$ satisfies 
$G(\cdot\, ; y)-\Phi(\cdot\, ; y)\in H^2_{loc}(\tilde{D})$ for all $y\in\tilde{D}$ and the open 
waveguide radiation condition; that is, has the decomposition in the form 
$$ G(\cdot\, ; y)\ =\ G_{rad}(\cdot\, ; y)\ +\ G_{prop}(\cdot\, ; y)\quad\mbox{in }\tilde{D}, $$
where $G_{prop}(\cdot\, ; y)$ is the propagating part. The radiating part $G_{rad}(\cdot\, ; y)$  includes the incoming wave and satisfies $G_{rad}(\cdot\, ; y)\in H^1(W_h\backslash\{B_\epsilon(y)\})$ for all $h>h_0$ and
some $\epsilon>0$ and also satisfies the Sommerfeld radiation condition. We shall prove the following properties of $G(x; y)$.

\begin{theorem}\label{sy}
\begin{itemize}
\item[(i)] The Green's function to the perturbed scattering problem satisfies $G(x; y)=G(y; x)$ 
for all $x,y\in\tilde{D}$, $x\neq y$.
\item[(ii)] In the unperturbed case (i.e., $\tilde{D}=D$), the propagating part $G_{prop}$ 
of $G$ takes the explicit form 
\begin{equation}\label{Gp}
G_{prop}(x; y)=2\pi i\sum_{j\in J} \biggl[ \psi_+(x_1)\sum_{\lambda_{\ell,j}>0} 
\frac{1}{\lambda_{\ell,j}}\hat{\phi}_{\ell,j}(x)\overline{\hat{\phi}_{\ell,j}(y)}
-\psi_-(x_1)\sum_{\lambda_{\ell,j}<0}\frac{1}{\lambda_{\ell,j}}\hat{\phi}_{\ell,j}(x) 
\overline{ \hat{\phi}_{\ell,j}(y)}\biggr]
\end{equation} 
for all $x,y\in D$, $x\neq y$.
\end{itemize}
\end{theorem}
\begin{remark}
In the perturbed case, the propagating part of $G$ can be decomposed into $G_{prop}=G_{prop}^{(0)}+G_{prop}^{(1)}$, where $G_{prop}^{(0)}$ 
represents the counterpart corresponding to the Green's function of the unperturbed problem 
taking the form \eqref{Gp}, while $G_{prop}^{(1)}$ denotes the propagating part caused by 
the defect.
\end{remark}
To prove Theorem \ref{sy} we need an auxiliary lemma. Below we take $R>\pi$ as a variable and suppose that $D_R:=\{x\in D:|x|<R\}$ is always a Lipschitz domain. Otherwise, we can slightly change the shape of the part $C_R$ to achieve this. In the remaining part of this paper we do not mention this any more.

\begin{lemma}\label{lem1}
Let 
\begin{eqnarray*}
u_{prop}(x) & = & \psi^+(x_1)\sum_{j\in J}\sum_{\lambda_{\ell,j}>0}a_{\ell,j}\,
\hat{\phi}_{\ell,j}(x)+ \psi^-(x_1)\sum_{j\in J}\sum_{\lambda_{\ell,j}<0}a_{\ell,j}\,
\hat{\phi}_{\ell,j}(x)\,, \\
v_{prop}(x) & = & \psi^+(x_1)\sum_{j\in J}\sum_{\lambda_{\ell,j}>0}b_{\ell,j}\,
\hat{\phi}_{\ell,j}(x)+ \psi^-(x_1)\sum_{j\in J}\sum_{\lambda_{\ell,j}<0}b_{\ell,j}\,
\hat{\phi}_{\ell,j}(x)\,,
\end{eqnarray*}
be the propagating parts of two solutions satisfying the open waveguide radiation condition. Then
$$ \int\limits_{C_R}\left[v_{prop}\,\frac{\partial u_{prop}}{\partial\nu}-
u_{prop}\,\frac{\partial v_{prop}}{\partial\nu}\right]ds\ \longrightarrow\ 0\quad
\mbox{as }R\to\infty\,. $$
\end{lemma}

\begin{remark}The path $C_R$ can be replaced by $\partial D_R$ because of the boundary condition on $\partial D$. Since 
$u_{prop},v_{prop}\in H^1(D_R)$ and $\Delta u_{prop},\Delta v_{prop}\in L^2(D_R)$ the integral 
over $\partial D_R$ is understood in the dual form of $\langle H^{-1/2}(\partial D_R),
H^{1/2}(\partial D_R)\rangle$. The integral over $C_R$ is understood in the dual form of 
$\langle H^{-1/2}(C_R),H^{1/2}_0(C_R)\rangle$ (see e.g., \cite{M10}).
\end{remark}
\textbf{Proof:} We first recall Green's second formula for any bounded Lipschitz 
domain $\Omega$. For $w,v\in H^1(\Omega)$ satisfying $\Delta w, \Delta v\in L^2(\Omega)$
we have
$$ \int\limits_\Omega\bigl[w\,\Delta v-v\,\Delta w\bigl]\, dx\ =\ 
\int\limits_{\partial\Omega}\bigl[w\,\partial_\nu v-v\,\partial_\nu w \bigl]\,ds $$
where the right hand sides are understood as the dual forms for $\partial_\nu v,
\partial_\nu w\in H^{-1/2}(\partial\Omega)$ and $w,v\in H^{1/2}(\partial\Omega)$.
Application of the Green's formula to $D_R$ yields
$$ \int\limits_{C_R}\left[v_{prop}\,\frac{\partial u_{prop}}{\partial\nu}-
u_{prop}\,\frac{\partial v_{prop}}{\partial\nu}\right]ds\ =\
\int\limits_{D_R}\bigl[\,v_{prop}\,(\Delta+k^2)u_{prop}-
u_{prop}\,(\Delta+k^2)v_{prop}\bigr]\,dx\,. $$
From the forms of $u_{prop}$ and $v_{prop}$ we conclude that
$$ \Delta u_{prop}+k^2u_{prop}\ =\
\sum_{j\in J}\sum_{\lambda_{\ell,j}>0}a_{\ell,j}\,(\Delta+k^2)[\psi^+\,
\hat{\phi}_{\ell,j}]\ +\ \sum_{j\in J}\sum_{\lambda_{\ell,j}<0}a_{\ell,j}\,
(\Delta+k^2)[\psi^-\,\hat{\phi}_{\ell,j}] $$
and the same for $v_{prop}$. Therefore, since $\psi^+(x_1)\psi^-(x_1)=0$ for all $x_1\in \R$,
\begin{eqnarray*}
& & \int\limits_{D_R}v_{prop}\,(\Delta+k^2)u_{prop}\,dx \\
& = & \sum_{j_1\in J}\sum_{\lambda_{\ell_1,j_1}>0}a_{\ell_1,j_1}
\sum_{j_2\in J}\sum_{\lambda_{\ell_2,j_2}>0}b_{\ell_2,j_2}
\int\limits_{D_R^+}[\psi^+\,\hat{\phi}_{\ell_2,j_2}]\,(\Delta+k^2)[\psi^+\,
\hat{\phi}_{\ell_1,j_1}]\,dx \\
& & + \sum_{j_1\in J}\sum_{\lambda_{\ell_1,j_1}<0}a_{\ell_1,j_1}
\sum_{j_2\in J}\sum_{\lambda_{\ell_2,j_2}<0}b_{\ell_2,j_2}
\int\limits_{D_R^-}[\psi^-\,\hat{\phi}_{\ell_2,j_2}]\,(\Delta+k^2)[\psi^-\,
\hat{\phi}_{\ell_1,j_1}]\,dx
\end{eqnarray*}
where $D_R^+=\{x\in D_R:\sigma-1<x_1<\sigma\}$ and $D_R^-=\{x\in D_R:-\sigma<x_1<-\sigma+1\}$.
Here we suppose that both  $D_R^+$ and  $D_R^-$ are Lipschitz domains by the choice of $\sigma>\max\{R, 2\pi\}+1$.
The analogous formula holds for $u_{prop}$ and $v_{prop}$ interchanged. Taking the 
difference and applying Green's theorem yields for the integral over $D_R^+$: 
\begin{eqnarray*}
& & \int\limits_{D_R^+}[\psi^+\,\hat{\phi}_{\ell_2,j_2}]\,(\Delta+k^2)[\psi^+\,
\hat{\phi}_{\ell_1,j_1}]-[\psi^+\,\hat{\phi}_{\ell_1,j_1}]\,(\Delta+k^2)[\psi^+\,
\hat{\phi}_{\ell_2,j_2}]\,dx \\ 
& = & \int\limits_{\partial D_R^+}\psi^+\left[\hat{\phi}_{\ell_2,j_2}
\frac{\partial[\psi^+\hat{\phi}_{\ell_1,j_1}]}{\partial\nu}-\hat{\phi}_{\ell_1,j_1}\,
\frac{\partial[\psi^+\hat{\phi}_{\ell_2,j_2}]}{\partial\nu}\right]ds \\
& = & \int\limits_{S_R}\psi^+\left[\hat{\phi}_{\ell_2,j_2}
\frac{\partial[\psi^+\hat{\phi}_{\ell_1,j_1}]}{\partial\nu}-\hat{\phi}_{\ell_1,j_1}\,
\frac{\partial[\psi^+\hat{\phi}_{\ell_2,j_2}]}{\partial\nu}\right]ds \\
& & +\ \int\limits_{\gamma_R}\left[\frac{\partial\hat{\phi}_{\ell_1,j_1}}{\partial x_1}\,
\hat{\phi}_{\ell_2,j_2}-\frac{\partial\hat{\phi}_{\ell_2,j_2}}{\partial x_1}\,
\hat{\phi}_{\ell_1,j_1}\right]ds
\end{eqnarray*}
where $S_R=\{x\in D:|x|=R,\ \sigma-1<x_1<\sigma\}$ and $\gamma_R=\{x\in D:x_1=\sigma,\ 
|x|<R\}$. We remark that, since $\hat{\phi}_{\ell,j}$ vanish on $\partial D$, the integral 
over $\gamma_R$ is understood  in the dual form of 
$\langle H^{-1/2}(\gamma_R),H^{1/2}_0(\gamma_R)\rangle$ (see e.g., \cite{M10}). The integral 
over $S_R$ tends to zero as $R\to\infty$ because of the exponential decay. The integral 
over $\gamma_R$ tends to
$$ \int\limits_{x_1=\sigma}\left[\frac{\partial\hat{\phi}_{\ell_1,j_1}}{\partial x_1}\,
\hat{\phi}_{\ell_2,j_2}-\frac{\partial\hat{\phi}_{\ell_2,j_2}}{\partial x_1}\,
\hat{\phi}_{\ell_1,j_1}\right]ds\ =\
\int\limits_{x_1=\sigma}\left[\frac{\partial\hat{\phi}_{\ell_1,j_1}}{\partial x_1}\,
\overline{\hat{\phi}_{\ell_2,-j_2}}-\frac{\partial\overline{\hat{\phi}_{\ell_2,-j_2}}}
{\partial x_1}\,\hat{\phi}_{\ell_1,j_1}\right]ds $$
because $\hat{\phi}_{\ell_2,-j_2}=\overline{\hat{\phi}_{\ell_2,j_2}}$. Now we conclude from 
$\lambda_{\ell_1,j_1}>0$ and $-\lambda_{\ell_2,-j_2}=\lambda_{\ell_2,j_2}>0$ that 
$j_1\not=-j_2$. Note that $\ell_1\in\{1,\ldots,m_{j_1}\}$ and $\ell_2\in
\{1,\ldots,m_{j_2}\}$ and $m_{j_2}=m_{-j_2}$. Therefore, the integral vanishes by the 
proof of \cite[Lemma~2.6]{K22}. \qed
\medlf
Below we carry out the proof of Theorem \ref{sy} by using Lemma \ref{lem1}. The symmetry of 
the Green's function will be used in the proof of Theorem \ref{Th-inverse-source} for our 
inverse problems. Write $B(x, \delta):=\{z\in \R^2: |z-x|<\delta\}$.
\medlf
\textbf{Proof of Theorem \ref{sy}:} (i) We fix $x,y\in\tilde{D}$ with $x\not=y$ and choose 
$\delta>0$ such that $B(x,\delta)\cup B(y,\delta)\subset\tilde{D}$ and 
$\overline{B(x,\delta)}\cap \overline{B(y,\delta)}=\emptyset$. Then we choose $R>0$ sufficiently large and set 
$$ D_{R,\delta}\ :=\ \bigl\{z\in\tilde{D}:|z|<R,\ |z-x|>\delta,\ |z-y|>\delta\bigr\}\,. $$
Using $\Delta_zG(z; x)+k^2G(z; x)=0$ and $\Delta_zG(z; y)+k^2G(z; y)=0$ in $D_{R,\delta}$ and 
application of Green's second formula in $D_{R,\delta}$ yields
\begin{eqnarray*}
0 & = & \int\limits_{D_{R,\delta}}\bigl[G(z; x)\,\Delta_zG(z; y)-G(z; y)\,\Delta_zG(z; x)
\bigr]\,dz \\
& = & \left(\int\limits_{|z|=R}-\int\limits_{|z-x|=\delta}-\int\limits_{|z-y|=\delta}\right)
\left[\frac{\partial G(z; y)}{\partial\nu(z)}\,G(z; x)-
\frac{\partial G(z; x)}{\partial\nu(z)}\,G(z; y)\right]ds(z)\,.
\end{eqnarray*}
Here we have used the vanishing of $G(\cdot\, ; y)$ and $G(\cdot\, ; x)$ on 
$\partial D_{R,\delta}\cap \tilde{\Gamma}$. Note that the normal direction at $C_R$ is supposed to direct into $\Sigma_R$ and that at $\partial B(x, \delta)$ or $\partial B(y, \delta)$ to direct into $D_{R,\delta}$.
We consider first the integral over $\partial B(x, \delta)$.
 For $|z-x|\leq\delta$ the terms $G(\cdot\, ; y)$ and 
$G(\cdot; x)-\Phi(\cdot; x)$ and their gradients are smooth. Therefore,
\begin{eqnarray*}
\lefteqn{ \int\limits_{|z-x|=\delta}\left[\frac{\partial G(z; y)}{\partial\nu(z)}\,G(z; x)-
\frac{\partial G(z; x)}{\partial\nu(z)}\,G(z; y)\right]ds(z) } \\
& = & \int\limits_{|z-x|=\delta}\left[\frac{\partial G(z; y)}{\partial\nu(z)}\,\Phi(z; x)-
\frac{\partial\Phi(z; x)}{\partial\nu(z)}\,G(z; y)\right]ds(z)\ +\ \mathcal{O}(\delta) \\
& = & G(x; y)\ +\ \mathcal{O}(\delta)
\end{eqnarray*}
where we applied Green's representation formula to $G(\cdot\, ; y)$ in the disk $B(x,\delta)$.
For $\delta\to 0$ we get
$$ \int\limits_{|z-x|=\delta}\left[\frac{\partial G(z; y)}{\partial\nu(z)}\,G(z; x)-
\frac{\partial G(z; x)}{\partial\nu(z)}\,G(z; y)\right]ds(z)\ \longrightarrow\ G(x; y)\,. $$
The integral over $\partial B(y, \delta)$ is treated in the same way, just by interchanging the roles 
of $x$ and $y$.
\smalf
It remains to show that the integral over $C_R:=\{z\in D:|z|=R\}$ tends to zero as $R$
tends to infinity. Substituting the decomposition $G(\cdot\, ; y)=G_{rad}(\cdot\, ; y)+
G_{prop}(\cdot\, ; y)$ into the integral yields that 
$$ \int\limits_{C_R}\left[\frac{\partial G(z; y)}{\partial\nu(z)}\,G(z; x)-
\frac{\partial G(z; x)}{\partial\nu(z)}\,G(z; y)\right]ds(z) $$
consists of four integrals. The integral
$$ \int\limits_{C_R}\left[\frac{\partial G_{prop}(z; y)}{\partial\nu(z)}\,G_{prop}(z; x)-
\frac{\partial G_{prop}(z; x)}{\partial\nu(z)}\,G_{prop}(z; y)\right]ds(z) $$
tends to zero by the previous lemma. For the other parts we note that the integral over 
$C_R\cap W_h$ tends to zero as $R\rightarrow\infty$, because, for example, 
\begin{eqnarray*}
& & \biggl|\int\limits_{C_R\cap W_h}\frac{\partial G_{rad}(z; y)}{\partial\nu(z)}
\,G_{prop}(z; x)\,ds(z)\biggr| \\
& \leq & c\,\biggl\Vert\frac{\partial G_{rad}(\cdot\, ; y)}{\partial\nu}
\biggr\Vert_{H^{-1/2}(C_R\cap W_h)}\,\Vert G_{prop}(\cdot,x)\Vert_{H^{1/2}(C_R\cap W_h)} \\ 
& \leq & c^\prime\,\Vert G_{rad}(\cdot\, ; y)\Vert_{H^1(V_R)}\,\Vert G_{prop}(\cdot,x)
\Vert_{H^1(V_R)}
\end{eqnarray*}
where $V_R:=\{z\in D\cap W_h:\bigl||z|-R\bigr|<1/2\}$,
$\Vert G_{rad}(\cdot\, ; y)\Vert_{H^1(U_R)}$ tends to zero and 
$\Vert G_{prop}(\cdot,x)\Vert_{H^1(V_R)}$ is bounded.

Hence it remains to consider the integral over $C_{R,h}:=\{z\in C_R:z_2>h\}$. Here we use 
Sommerfeld's radiation condition for $G_{rad}(\cdot,x)$ and $G_{rad}(\cdot\, ; y)$ which yields 
that 
$$ \int\limits_{C_{R,h}}\left[\frac{\partial G_{rad}(z; y)}{\partial\nu(z)}\,G_{rad}(z; x)-
\frac{\partial G_{rad}(z; x)}{\partial\nu(z)}\,G_{rad}(z; y)\right]ds(z) $$
tends to zero. Using the estimates 
\ben
|G_{rad}(z; x)|+|\nabla_z G_{rad}(z; x)|\leq\frac{c}{\sqrt{|z|}},\quad 
|G_{prop}(z; x)|+|\nabla_z G_{prop}(z; x)|\leq c\,e^{-\sigma z_2} 
\enn
for some $c>0$ 
and the same for $y$ replacing $x$, we obtain
\begin{eqnarray}\nonumber
& & \int\limits_{C_{R,h}}\left|\frac{\partial G_{prop}(z; y)}{\partial\nu(z)}\,G_{rad}(z; x)-
\frac{\partial G_{rad}(z; x)}{\partial\nu(z)}\,G_{prop}(z; y)\right|ds(z) \\ \nonumber
& & +\ \int\limits_{C_{R,h}}\left|\frac{\partial G_{rad}(z; y)}{\partial\nu(z)}\,G_{prop}(z; x)-
\frac{\partial G_{prop}(z; x)}{\partial\nu(z)}\,G_{rad}(z; y)\right|ds(z)\\ \label{R}
& \leq&\ 
c\,\frac{R}{\sqrt{R}}\int\limits_0^\pi e^{-R\sigma\sin t}\,dt. 
\end{eqnarray}
Using $\sin t\geq \frac{2}{\pi}t$ on $[0,\pi/2]$, one deduces that the right hand side of \eqref{R} 
tends to zero as $R\rightarrow\infty$, because $$\int\limits_0^\pi e^{-R\sigma\sin t}\,dt=2\int\limits_0^{\pi/2}
e^{-R\sigma\sin t}\,dt\leq 2\int\limits_0^{\pi/2}
e^{-2R\sigma t/\pi}\,dt
=\frac{\pi}{R\sigma}(1-e^{-R\sigma}).$$ 
\medlf
(ii) For fixed $y\in D$, we choose $\epsilon>0$ less than the distance between 
$y$ and $\Gamma$. Introduce a cut-off function $\chi\in C_0^\infty(\real^2)$ with $\chi(x)=1$ 
for $|x-y|<\epsilon/2$ and $\chi(x)=0$ for $|x-y|\geq\epsilon$. Then 
$v:=G(\cdot\, ; y)-\chi\,\Phi(\cdot\, ; y)\in H^1_{loc}(D)$ coincides with $G(\cdot\, ; y)$ for 
$|x-y|\geq\epsilon$ and satisfies $\Delta v+k^2v=-g_y$ in $D$ and $v=0$ on $\Gamma$ where
$$ g_y\ :=\ \Delta\chi\,\Phi(\cdot\, ; y)+2\nabla\chi\cdot\nabla\Phi(\cdot\, ; y)\in L^2(D) $$
has compact support. If $\sigma$ in the definition of the radiation condition is chosen to 
be larger than $|y_1|+\epsilon$ then, by \cite[Theorem 3.5]{HK22},
\begin{eqnarray*}
G_{prop}(x; y) & = & v_{prop}(x)\ =\ \sum_{j\in J}\biggl[\psi_+(x_1)\sum_{\lambda_{\ell,j}>0}
\hat{\phi}_{\ell,j}(x)\,a_{\ell,j}(y)+\psi_-(x_1) \sum_{\lambda_{\ell,j}<0}
\hat{\phi}_{\ell,j}(x)\,a_{\ell,j}(y)\biggr]
\end{eqnarray*} 
for $|x_1|\geq\sigma $ where the coefficients $a_{\ell,j}(y)$ are given by
$$ a_{\ell,j}(y)\ =\ \frac{2\pi i}{|\lambda_{\ell,j}|}\int\limits_{\epsilon/2<|x-y|<\epsilon} 
g_y(x)\,\overline{\hat{\phi}}_{\ell,j}(x)\,dx\,. $$
To calculate $a_{\ell,j}(y)$, we rewrite $g_y$ as $g_y(x)=(\Delta+k^2)(\chi(x)\Phi(x ; y))$ for $x\neq y$. 
Consequently, application of Green's formula yields
\begin{eqnarray}\nonumber
a_{\ell,j}(y) & = & \frac{2\pi i}{|\lambda_{\ell,j}|}\int\limits_{\epsilon/2<|x-y|<\epsilon} 
g_y(x)\,\overline{\hat{\phi}_{\ell,j}(x)}\,dx \\ \nonumber
& = & \frac{2\pi i}{|\lambda_{\ell,j}|}\int\limits_{|x-y|=\epsilon/2}\biggl[\frac{\partial
\overline{\hat{\phi}_{\ell,j}(x)}}{\partial\nu(x)}\,\Phi(x; y)-\frac{\partial\Phi(x; y)}
{\partial\nu(x)}\,\overline{\hat{\phi}_{\ell,j}(x)}\biggr]\,ds(x) \\ \label{alj}
& = & \frac{2\pi i}{|\lambda_{\ell,j}|}\,\overline{\hat{\phi}_{\ell,j}(y)}\,,
\end{eqnarray}
where we have used the fact that $\chi(x)=1$ on $|x-y|=\epsilon/2$ and that $\chi(x)=0$ on $|x-y|=\epsilon$.  
\qed
\smalf
From the proof of Theorem \ref{sy}, we conclude that
\begin{corollary}
Let $u\in H^1_{loc, 0}(\tilde{D})$ be an open wave-guide radiating solution to the Helmholtz equation $(\Delta+k^2)u=0$ in $\Sigma_R$, where $\tilde{D}_R$ is supposed to be a Lipschitz domain.  We have the representation
\be\label{eq:u}
u(x)=\int\limits_{C_R}
\left[\frac{\partial u(z)}{\partial\nu(z)}\,G(z; x)-
\frac{\partial G(z; x)}{\partial\nu(z)}\,u(z)\right]ds(z),\quad x\in \Sigma_R.
\en
\end{corollary}
\begin{proof}
We fix $x\in \Sigma_R$ and choose $R'>R$ such that $x\in \tilde{D}_{R'}$ where $\tilde{D}_{R'}$ is a Lipschitz domain. Application of the Green's representation formula yields
\be\label{exu}
u(x)=\left(\int\limits_{C_R}-\int\limits_{C_{R'}}\right)
\left[\frac{\partial u(z)}{\partial\nu(z)}\,G(z; x)-
\frac{\partial G(z; x)}{\partial\nu(z)}\,u(z)\right]ds(z),\quad x\in \Sigma_R.
\en
By the proof of Theorem \ref{sy}, the integral over $C_{R'}$ tends to zero as $R'$ tends to infinity, which together with \eqref{exu} finishes the proof of \eqref{eq:u}.
\end{proof}

\section{Scattering of Plane Waves at a Propagative Wave Number}\label{sec3-1}

As shown in Proposition \ref{wpw}, uniqueness and existence of a weak solution for an incoming 
plane wave are guaranteed under the open waveguide radiation condition, provided 
$k\sin\theta$ is not a propagative wave number. If $k\sin\theta=\hat{\alpha}_j$ is a 
propagative wave number, there exists still a $\hat{\alpha}_j$-quasi-periodic solution 
$u_0\in H^1_{loc}(D)$ of the unperturbed problem (see Lemma~\ref{lem3.4} (i) below). 
However, uniqueness fails and the general solution takes the form
\begin{equation}\label{unpert-1}
u\ =\ u_0\ +\ \sum_{\ell=1}^{m_j}c_\ell\,\hat{\phi}_{\ell, j}\quad\mbox{in }D
\end{equation} 
where $\hat{\phi}_{\ell,j}\in X_j$ (see \eqref{X_j} and \eqref{evp}) and $c_\ell\in\mathbb{C}$
are arbitrary. A general solution to the locally perturbed scattering problem is described in 
\cite[Corollary 4.10]{HK22} when $k\sin\theta$ is a propagative wave number. The purpose of 
this section is to propose an additional constraint on solutions of the unperturbed problem 
to fix the coefficients $c_\ell$ in \eqref{unpert-1}, so that the forward problem is always 
uniquely solvable. We shall employ two methods: The Limiting Absorption Principle for approximating the wave number with a positive imaginary part in Subsection \ref{LAP} and the method of approximating the plane wave by point source waves in Subsection \ref{APS}.  

\subsection{The Limiting Absorption Principle and Singular Perturbation Arguments}
\label{LAP}

We first consider the scattering problem in a periodic domain $D$ without defects.
Let $u^{in}(x)=e^{ik(x_1\sin\theta-x_2\cos\theta)}$ with $\theta\in(-\pi/2, \pi/2)$ be the 
incident plane wave. Set 
\ben
\alpha\ :=\ k\sin\theta\quad\mbox{and}\quad \beta_n\ :=\ \sqrt{k^2-(n+\alpha)^2}\,, \quad n\in \N,
\enn
where the square root is chosen such that $\Im \sqrt{t}\geq 0$ for $t\leq0$. Obviously, $u^{in}(x)=e^{i\alpha x_1-i\beta_0x_2}$.
Since the wave number $k>0$ is fixed, we omit the dependence on $k$ for simplicity. We look
for an $\alpha$-quasi-periodic total field $u\in H^1_{\alpha,0}(Q_h):=\{u\in H^1_\alpha(Q_h): u=0\;\mbox{on}\;\partial Q_h\cap \Gamma\}$ for all 
$h>h_0$ such that $u^{sc}=u-u^{in}$ satisfies the upward $\alpha$-quasiperiodic Rayleigh expansion \eqref{exc:b}.

Introduce the $\alpha$-quasi-periodic and periodic, respectively, Sobolev spaces on $\G_h$ 
with $h>h_0$ by
\ben
H_\alpha^{1/2}(\G_h) & := & \{f\in H^{1/2}(\G_h):e^{-i\alpha x_1}f(x_1,h)
\mbox{ is $2\pi$-periodic in $x_1$}\}\,, \\
H_{per}^{1/2}(\G_h) & := & \{f\in H^{1/2}(\G_h):f(x_1,h)\mbox{ is $2\pi$-periodic in $x_1$}\}\,.
\enn
Define the periodic and quasi-periodic, respectively, Dirichlet-to-Neumann maps on the 
artificial boundary $\Gamma_h$ by 
\begin{eqnarray}
\hspace*{1cm}(T_kf)(x_1,h) & := & \sum_{n\in\Z}i\beta_n\,f_n\,e^{inx_1},\quad f(x_1,h)=
\sum_{n\in\Z}f_n\,e^{inx_1}\in H_{per}^{1/2}(\G_h)\,, \label{dtn} \\ 
\hspace*{1cm}(\tilde{T}_{k}\tilde{f})(x_1,h) & := & \sum_{n\in\Z}i\beta_n\,\tilde{f}_n\,
e^{i(n+\alpha)x_1},\quad \tilde{f}(x_1,h)=\sum_{n\in\Z}\tilde{f}_n 
e^{i(n+\alpha)x_1}\in H_\alpha^{1/2}(\G_h)\,. \label{dtn2}
\end{eqnarray}
It is well known that  $T_k: H_{per}^{1/2}(\Gamma_h)\rightarrow H_{per}^{-1/2}(\Gamma_h)$ 
and $\tilde{T}_k: H_{\alpha}^{1/2}(\Gamma_h)\rightarrow H_{\alpha}^{-1/2}(\Gamma_h) $ are 
bounded linear operators. The following variational formulation for $u\in H^1_{\alpha,0}(Q_h)$ can be easily derived:
\be\label{var-u}
\tilde{a}_k(u, \phi)\ =\ -2ik\cos\theta\, e^{-ikh\cos\theta}\int\limits_{\G_h} 
e^{i\alpha x_1}\,\overline{\phi}\,ds\quad\mbox{for all }\phi\in H^1_{\alpha,0}(Q_h)\,,
\en
where 
$$ \tilde{a}_k(u,\phi)\ :=\ \int\limits_{Q_h}\bigl[\nabla u\cdot\overline{\nabla\phi}-
k^2 u\,\overline{\phi}\bigr]\,dx\ -\ \int\limits_{\G_h}\tilde{T}_k u\,\overline{\phi}\,ds\,,
\quad u,\phi\in H^1_{\alpha,0}(Q_h)\,. $$
Defining
$v:=e^{-i\alpha x_1}u\in H^1_{per,0}(Q_h):=\{v\in H^1_{per}(Q_h), v=0\mbox{ on }\Gamma\cap \partial Q_h\}$  
and $\psi:=e^{-i\alpha x_1}\phi$ for $\phi\in H^1_{\alpha,0}(Q_h)$, we get the periodic
form 
\be\label{eq8}
a_k(v,\psi)\ =\ -2ik\cos\theta\,e^{-ikh\cos\theta}\int\limits_0^{2\pi}\overline{\psi(x_1,h)}
\,dx_1\quad\mbox{for all }\quad\psi\in H^1_{per,0}(Q_h)
\en
where now 
$$ a_k(v,\psi)\ :=\ \int\limits_{Q_h} \biggl[\nabla v\cdot\overline{\nabla\psi}-
2i\alpha\,\frac{\partial v}{\partial x_1}\,\overline{\psi}-(k^2-\alpha^2)\,v\,
\overline{\psi}\biggr]\,dx\ -\ \int\limits_{\G_h} T_k v\,\overline{\psi}\,ds$$
for $v,\psi\in H^1_{per,0}(Q_h)$. 
We equip $H^1_{per,0}(Q_h)$ with the inner product
\be\label{InPr}
\langle v,\psi\rangle\ :=\ \int\limits_{Q_h}\nabla v\cdot\nabla\overline{\psi}\,dx\ +\
2\pi\sum_{n\in\Z}(1+n^2)^{1/2}v_n\,\overline{\psi}_n\,,\quad\psi,v\in H^1_{per,0}(Q_h)\,,
\en
where $\psi_n$ and $v_n$ denote the Fourier coefficients of $\psi(x_1,h)$ and $v(x_1, h)$, 
respectively. By the representation theorem of Riesz, there exist $f^{(k)}\in 
H^1_{per,0}(Q_h)$ and a linear bounded operator $L_k$ from $H^1_{per,0}(Q_h)$ into itself with 
\begin{eqnarray}
\langle f^{(k)},\psi\rangle & = & -2ik\cos\theta\,e^{-ikh\cos\theta}\int\limits_0^{2\pi}
\overline{\psi(x_1,h)}\,dx_1\,, \label{fk} \\
\langle L_kv,\psi\rangle & = & a_k(v,\psi) \label{Lk} \\ 
& = & \int\limits_{Q_h}\biggl[\nabla v\cdot\nabla\overline{\psi}-2i k\sin\theta\,
\frac{\partial v}{\partial x_1}\,\overline{\psi}-k^2\cos^2\theta\,v\,\overline{\psi}\biggr]\,
dx\ -2\pi\sum_{n\in \N} i\beta_n v_n\overline{\psi}_n  \nonumber
\end{eqnarray}
for all $\psi,v\in H^1_{per,0}(Q_h)$. Then the operator equation \eqref{eq8} can be rewritten 
as
\be\label{eq7}
L_kv\ =\ f^{(k)}\quad\mbox{in}\quad H^1_{per,0}(Q_h)\,.
\en 
Using the compact embedding of $H^1_{per,0}(Q_h)\rightarrow L^2(Q_h)$, one can show that the operator $K_k:=I-L_k$, given by
\ben
\langle K_k\,v,\psi\rangle
& := & \int\limits_{Q_h}\biggl[2i k\sin\theta\,
\frac{\partial v}{\partial x_1}\,\overline{\psi}+k^2\cos^2\theta\,v\,\overline{\psi}\biggr]\,
dx\ +2\pi\sum_{n\in \N}(\sqrt{1+n^2}+ i\beta_n) v_n\overline{\psi}_n,  
\enn
for all $\psi,v\in H^1_{per,0}(Q_h)$,
 is compact as an operator from $H^1_{per,0}(Q_h)$ into 
itself. 
\medlf
Below we collect some properties of the operator $L_k$.

\begin{lemma}\label{lem3.4} Suppose that $\alpha+n\neq \pm k$ for any $n\in \Z$, that is, 
$\alpha$ is not a cut-off value. 
\begin{itemize}
\item[(i)] For $k>0$, the equation \eqref{eq7} admits at least one solution 
$v\in H^1_{per,0}(Q_h)$. The null space $\mathcal{N}:=\mathcal{N}(L_k)=\mathcal{N}(L_k^*)$ 
is finite dimensional and consists of surface wave modes only, i.e., 
\be\label{eq12}
v(x)\ =\ \sum_{n\in\Z:|n+\alpha|>k}v_n\,e^{inx_1-|\beta_n|(x_2-h)}\,,\quad x_2>h\,.
\en
\item[(ii)] The Riesz number of $L_k$ is one, that is, $\mathcal{N}(L_k)=\mathcal{N}(L_k^2)$. 
Moreover, there holds the orthogonal decomposition $H^1_{per,0}(Q_h)=\mathcal{N}(L_k)\oplus
\mathcal{R}(L_k)$. Here $\mathcal{R}(L_k)$ denotes the range of the operator $L_k$.
\item[(iii)] If $\Im k>0$, there is a unique solution to \eqref{eq7} in $H^1_{per,0}(Q_h)$ for any $h>h_0$.
\end{itemize}
\end{lemma}
\begin{proof}  (i) The form of $v\in\mathcal{N}$ given by \eqref{eq12} can be derived by 
setting $\phi=u=v e^{i\alpha x_1}$ in the homogeneous form of equation \eqref{var-u}, taking the 
imaginary part and using the definition of $\tilde{T}_k$. The adjoint operator $L_k^\ast$ 
of $L_k$ is defined by
\begin{eqnarray*}
\langle L_k^\ast v,\psi \rangle & = & \langle v, L_k\psi\rangle\ =\ 
\overline{\langle L_k\psi,v\rangle}\ =\ \overline{a_k(\psi,v)} \\
& = & \int\limits_{Q_h}\left[\nabla u\cdot\overline{\nabla \phi}-k^2 u\overline{\phi}
\right]\,dx\ +\ \sum_{n\in\Z} i\overline{\beta_n}\,u_n\,\overline{\phi}_n
\end{eqnarray*}
for all $v,\psi\in H^1_{per,0}(Q_h) $, where $u=e^{i\alpha x_1}v$ and 
$\phi=e^{i\alpha x_1}\psi$. From this we conclude that $\mathcal{N}(L_k)=
\mathcal{N}(L_k^\ast)$. The existence of $v\in H^1_{per,0}(Q_h)$ follows from the fact that  
$\langle f^{(k)}, \psi\rangle=0$ for all $\psi\in \mathcal{N}$ and the Fredholm alternative. 
The space $\mathcal{N}$ is finite dimensional, because $K_k$ is compact. 
\smalf
(ii) It is obvious that $\mathcal{N}(L_{k})\subset\mathcal{N}(L_k^2)$. To prove the reverse 
direction, we assume $L_k^2 w=0$ for some $w\in H^1_{per,0}(Q_h)$  and set 
$v=L_k w\in\mathcal{R}(L_k)$. Since $v\in\mathcal{N}(L_k)=\mathcal{N}(L_k^\ast)$, we obtain 
\ben
\Vert v\Vert^2\ =\ \langle v,v\rangle\ =\ \langle v,L_kw\rangle\ =\ 
\langle L_k^\ast v,w\rangle\ =\ 0, 
\enn
which proves $\mathcal{N}(L_{k}^2)\subset\mathcal{N}(L_k)$ and thus the coincidence 
$\mathcal{N}(L_{k}^2)=\mathcal{N}(L_k)$. This also implies $\mathcal{N}\cap\mathcal{R}=
\emptyset$ and hence $H^1_{per,0}(Q_h)=\mathcal{N}\oplus\mathcal{R}$. The orthogonality 
between $\mathcal{N}$ and $\mathcal{R}$ follows from the relation 
$\mathcal{N}(L_k)=\mathcal{N}(L_k^\ast)$.
\smalf
(iii) We shall apply the uniqueness result of \cite{CR1995} to the proof of the third assertion.
By Fredholm alternative, it suffices to prove uniqueness. Let $k\in \C$ with $\Im k>0$ and 
set $\alpha=k\sin\theta\in \C$. Assuming $L_kv=0$ for some $v\in H^1_{per,0}(Q_h)$, we need 
to prove $v\equiv 0$. It then follows that $a_k(v,\psi)=0$ for all $\psi\in H^1_{per,0}(Q_h)$, 
which implies that $v$ satisfies the elliptic equation 
$(\Delta+2i\alpha \partial_1 + k^2-\alpha^2) v=0$ in $D$, the Dirichlet boundary condition 
$v=0$ on $\Gamma$ together with the periodic expansion 
$v=\sum_{n\in \Z} v_n e^{in x_1+i\beta_n x_2}$ in $x_2>h$. 

We claim that $\Im\beta_n>0$ for all 
$n\in\ganz$, if $\Im k>0$. Recall that 
 the square root function $z\mapsto\sqrt{z}$ was chosen to be 
holomorphic in the region $\{z\in\cmplx: z\notin i\real_{\leq 0}\}$. In particular we have that 
$\Im\sqrt{z}>0$ for all $z\in C:=\{z\in\cmplx:\Re z<0\mbox{ or }\Im z>0\}$. It is easily seen that 
$k^2-(n+k\sin\theta)^2\in C$ for all $k\in\cmplx$ with $\Im k>0$ and $\Re k>0$ and all $n\in\ganz$
provided that $|\sin\theta|<1$. Indeed, if $n$ is such that $|n+\Re k\sin\theta|\leq\Re k$, then
\ben
\Im[k^2-(n+k\sin\theta)^2]&=&2\Im k\, [\Re k-(n+\Re k\sin\theta)\sin\theta]\\
&\geq& 
2\Im k\Re k\,(1-|\sin\theta|) 
>0.
\enn If $n$ is such that $|n+\Re k\sin\theta|>\Re k$, then
\ben
\Re[k^2-(n+k\sin\theta)^2]&=&(\Re k)^2-(\Im k)^2-(n+\Re k\sin\theta)^2+(\Im k)^2\sin^2\theta\\
&=&
(\Re k)^2-(n+\Re k\sin\theta)^2-(\Im k)^2\cos^2\theta
<0.
\enn This proves $\Im\beta_n>0$ for all 
$n\in\ganz$. 

Setting $u^s=v\, e^{i\alpha x_1}$, we deduce that $u^s$ fulfills the homogeneous boundary value
problem of the Helmholtz equation
$$ \Delta u^s+k^2 u^s=0\quad\mbox{in}\quad D,\qquad u^s=0\quad\mbox{on}\quad \partial D, $$
and the quasi-periodic Rayleigh expansion condition \eqref{exc:b} in $x_2>h$. 
Since $\Im\beta_n>0$ the function $u^s$ decays exponentially to zero as $x_2\rightarrow\infty$ 
in $D$. By elliptic boundary regularity in Lipschitz domains with the zero boundary condition, we have $u^s\in C^{0,\gamma}(\overline{Q_h})$ for any $h>h_0$, where the H\"older exponent $\gamma>0$ depends on the Lipschitz constant of $\Gamma$; we refer to \cite[Chapter 4, Theorem 4.3]{Gri16} for the proof valid for $C^1$-smooth boundaries, which is also applicable to Lipschitz curves in $\R^2$. The boundary behavior of the Dirichlet boundary value problem of elliptic equations in a non-smooth domain can be further found in \cite[Chapter 7]{Kry96}. Hence, this together with the interior regularity gives $u^s\in C^2(D)\cap C(\overline{D})$ and $u^s$ 
must be bounded in the infinite strip $Q_\infty:=\{x\in D: 0<x_1<2\pi\}$. On 
the other hand,  the $\alpha$-quasi-periodicity of $u^s$ gives
$$ u^s(x_1+2n\pi,x_2)=e^{i2\pi nk\sin\theta}u^s(x_1,x_2)\quad\mbox{for all }
x\in D.$$
This in combination with the boundedness of $\Vert u\Vert|_{L^\infty(Q_\infty)}$ yields
the growth condition
$$ |u^s(x)|\leq C\, e^{\Im k\,|\sin\theta|\,|x|}\,,\quad x\in D $$
for some $C>0$. Now, applying \cite[Theorem 3.1]{CR1995} we get $u^s\equiv 0$ in $D$ and thus 
$v=0$ in $H^1_{per,0}(Q_h)$ for all $h>h_0$.

\end{proof}
\medlf
Now we suppose that $k\sin\theta=\hat{\alpha}_j$ for some $j\in J$ is a propagative wave number, which implies 
$\mathcal{N}\neq\emptyset$. Replacing $k$ by $k+i\epsilon$ with $\epsilon>0$, we consider 
the perturbed operator equation
\ben
L(\epsilon)\,v_\epsilon\ =\ f(\epsilon)\quad\mbox{in}\quad H^1_{per,0}(Q_h)
\enn 
with
$L(\epsilon)=L_{k+i\epsilon}$ and $f(\epsilon)=f^{(k+i\epsilon)}$.  We want to study the 
convergence and limit of $v(\epsilon)$ as $\epsilon\rightarrow 0^+$ by applying the following 
singular perturbation result from \cite[Theorem 2.7 and Remark 2.8]{K22p}.
\begin{lemma}\label{sipe}
Let $I=(0, \epsilon_0)$ for some small $\epsilon_0>0$. Let $K(\epsilon)$ be compact operators from 
some Hilbert space $X$ into itself and $f(\epsilon)\in \mathcal{R}(L(\epsilon))$ for all $\epsilon\in[0, \epsilon_0)$
where $L(\epsilon):=I-K(\epsilon)$. Furthermore, let $L(\epsilon)$ be 
one-to-one (thus invertible) for all $\epsilon\in I$ and let $L(0)=I-K(0)$ have Riesz number 
one. Let $P:X\rightarrow\mathcal{N}:=N(L(0))$ be the projection onto the nullspace of $L(0)$ 
along the direct decomposition $X=\mathcal{N}\oplus\mathcal{R}$ where 
$\mathcal{R}=\mathcal{R}(L(0))$. Finally, let $f(\epsilon)$ and $K(\epsilon)$ be continuously 
differentiable functions in $\epsilon\in [0,\epsilon_0)$ and let $PL^\prime(0)|_{\mathcal{N}}$ 
be an isomorphism from $\mathcal{N}$ onto itself where $L^\prime(0)$ denotes the one-sided 
derivative of $L(\epsilon)$ at $\epsilon=0^+$.

Then the mapping $\epsilon\mapsto v(\epsilon):=[L(\epsilon)]^{-1}f(\epsilon)$ has a 
continuous  extension to $[0,\epsilon_0)$ into $X$. The limit 
$v(0)=\lim_{\epsilon\rightarrow 0+}v(\epsilon)$ is the unique solution of the system
\be\label{elimit}
L(0)\,v(0)\ =\ f(0)\,,\quad [PL^\prime(0)]\,v(0)\ =\ Pf^\prime(0)\,,
\en where $f^\prime(0)$ denotes the right hand derivative of $f(\epsilon)$ at $\epsilon=0$. Moreover, 
there exist $\delta\in(0,\epsilon_0)$ and $c>0$ such that
\ben
\Vert v(\epsilon_1)\Vert_X\ \leq\ c\,\bigl[\sup_{\epsilon\in[0,\delta]}
\Vert f(\epsilon)\Vert_X+\sup_{\epsilon\in[0,\delta]}
\Vert f^\prime(\epsilon)\Vert_X\bigr]\quad\mbox{for all }\epsilon_1\in[0,\delta]\,.
\enn
\end{lemma}
The original version of Lemma \ref{sipe} can be found in \cite[Theorem 1.32 , Section 1.4]{CK13}. A more direct proof is presented in \cite{K22p} with the characterization of the equation \eqref{elimit} of the limiting solution. 
To apply Lemma \ref{sipe} to the operator equation \eqref{eq7}, we set 
$X=H^1_{per,0}(Q_h)$ and denote by $P: H^1_{per,0}(Q_h)\rightarrow\mathcal{N}$ the projection 
operator. For $\psi\in H^1_{per,0}(Q_h)$, it follows from the definitions of $f^{(k)}$ 
and $L_k$ that (e.g., \eqref{fk} and \eqref{Lk})
\begin{eqnarray}
\langle f(\epsilon),\psi\rangle & = & -2i(k+i\epsilon)\cos\theta\,
e^{-i(k+i\epsilon)h\cos\theta}\int\limits_0^{2\pi}\overline{\psi(x_1,h)}\,dx_1,\nonumber \\
\langle L(\epsilon)v,\psi\rangle & = & \int\limits_{Q_h}\biggl[
\nabla v\cdot\nabla\overline{\psi}-2i(k+i\epsilon)\sin\theta\,\frac{\partial v}{\partial x_1}\, 
\overline{\psi}-(k+i\epsilon)^2\cos^2\theta v\overline{\psi}\biggr]\,dx \nonumber \\ 
\label{eq10}
& & -\ 2\pi \sum_{n\in \Z} i\sqrt{(k+i\epsilon)^2-[n+(k+i\epsilon)\sin\theta]^2}\,v_n\,
\overline{\psi_n}
\end{eqnarray}
where $v(x_1,h)=\sum_{n\in \Z}v_ne^{in x_1}\in X$ and $\psi(x_1,h)=\sum_{n\in \Z}\psi_ne^{in x_1}\in X$.
From the above expressions we observe that $f$ and $L$ are differentiable with respect to $\epsilon$ in a 
neighborhood of $0$ provided $k\sin\theta$ is not a cut-off value. By Lemma \ref{lem3.4} 
(iii),  $L(\epsilon)$ is invertible  and thus $f(\epsilon)\in\mathcal{R}(L(\epsilon))$ for all 
$\epsilon>0$. On the other hand, we have $f(0)\in \mathcal{R}$, because by Lemma \ref{lem3.4} (i) and (ii), $f(0)$ is orthogonal to $\mathcal{N}$ and $X$ admits the orthogonal decomposition $X=\mathcal{N}\oplus\mathcal{R}$.

 Since the null space $\mathcal{N}$ consists of evanescent wave modes only and 
$(I-P)h$ is orthogonal to $\mathcal{N}$ for any $h\in H^1_{per,0}(Q_h)$,  it holds that
\ben
\langle Pf(0),\psi\rangle & = & \langle f(0),\psi\rangle\ =\ -2ik\cos\theta\,
e^{-ikh \cos\theta}\int\limits_0^{2\pi}\overline{\psi(x_1, h)}\, dx_1\ =\ 0\,, \\
\langle Pf^\prime(0),\psi\rangle & = &\langle f^\prime(0),\psi\rangle= 2\cos\theta(1-ikh\cos\theta)\,e^{-ikh\cos\theta}
\int\limits_0^{2\pi}\overline{\psi(x_1, h)}\, dx_1\ =\ 0
\enn
for all $\psi\in\mathcal{N}$. This implies that $f(0),f^\prime(0)\in\mathcal{R}$ and 
$Pf(0)=Pf^\prime(0)=0$. On the other hand, simple calculations show for 
$v,\psi\in \mathcal{N}$ that
\be\nonumber
\langle PL^\prime(0)v,\psi\rangle  = \langle L^\prime(0)v,\psi\rangle&=& \int\limits_{Q_h}\biggl[2\sin\theta\,\frac{\partial v}
{\partial x_1}\,\overline{\psi}-2ik\cos^2\theta v\overline{\psi}\biggr]\,dx \\ \label{PK0}
& & +2i\pi\sum_{n\in \Z: |n+\alpha|>k}\frac{(n+\alpha)\sin\theta-k}{\sqrt{(n+\alpha)^2-k^2}}\,
v_n\,\overline{\psi}_n\,.
\en
It remains to justify the one-to-one property of the mapping $PL'(0)|_{\mathcal{N}}$ from 
$\mathcal{N}$ onto itself, which is given by the lemma below.

\begin{lemma}\label{lem3.6}
$PL^\prime(0)$ is one-to-one on $\cN$.
\end{lemma}
\textbf{Proof:} First we show that
\begin{equation}\label{PK} 
\langle PL^\prime(0)v,\psi\rangle\ =\ 2\int\limits_{Q_\infty}\biggl[
\sin\theta\,\frac{\partial v}{\partial x_1}\,\overline{\psi}- 
ik\cos^2\theta\,v\,\overline{\psi}\biggr]\,dx 
\end{equation}
for all $v,\psi\in\cN$ where $v$ and $\psi$ are extended into $Q_\infty\setminus Q_h$ by
\begin{eqnarray*}
v(x) & = & \sum_{n\in \Z: |n+\alpha|>k}v_n\,e^{-\sqrt{(n+\alpha)^2-k^2}(x_2-h)+i nx_1}\,,\quad 
x_2>h\,, \\
\psi(x) & = & \sum_{n\in \Z: |n+\alpha|>k}\psi_n\,e^{-\sqrt{(n+\alpha)^2-k^2}(x_2-h)+i nx_1},\,\quad 
x_2>h\,.
\end{eqnarray*}
Indeed, we compute
\begin{eqnarray*}
& & 2\int\limits_{Q_\infty\setminus Q_h}\biggl[\sin\theta\,\frac{\partial v}{\partial x_1}\,
\overline{\psi}-ik\cos^2\theta\,v\,\overline{\psi}\biggr]\,dx \\
& = & 4\pi i\int\limits_h^\infty\sum_{n\in \Z: \,|n+\alpha|>k}v_n\,\overline{\psi_n}\,(\sin\theta\,n-
k\cos^2\theta)\,e^{-2\sqrt{(n+\alpha)^2-k^2}(x_2-h)}dx_2 \\
& = & 2\pi i\sum_{n\in \Z: \,|n+\alpha|>k}v_n\,\overline{\psi_n}\,\frac{\sin\theta\,n-k\cos^2\theta}
{\sqrt{(n+\alpha)^2-k^2}} \\
& = & 2\pi i\sum_{n\in \Z: \,|n+\alpha|>k}v_n\,\overline{\psi_n}\,\frac{(n+\alpha)\sin\theta-k}
{\sqrt{(n+\alpha)^2-k^2}}\,.
\end{eqnarray*}
This in combination with \eqref{PK0} yields \eqref{PK}.
\medlf
Assume now that $\langle PL^\prime(0)v,\cdot\rangle$ vanishes identically on $\cN$ for some 
$v\in\cN$. Then $\langle PL^\prime(0)v,v\rangle=0$ and thus
\begin{eqnarray}\label{re1}
\sin\theta\int\limits_{Q_\infty}\frac{\partial v}{\partial x_1}\,\overline{v}\,dx\ =\ 
ik\cos^2\theta\int\limits_{Q_\infty}|v|^2dx\,. 
\end{eqnarray}
We substitute $v(x)=e^{-ik\sin\theta x_1}u(x)$ again and have
$$ \sin\theta\int\limits_{Q_\infty}\biggl[-ik\sin\theta\,|u|^2+\frac{\partial u}{\partial x_1}\,
\overline{u}\biggr]\,dx\ =\ ik\cos^2\theta\int\limits_{Q_\infty}|u|^2dx\,; $$
that is,
$$ \sin\theta\Im\int\limits_{Q_\infty}\frac{\partial u}{\partial x_1}\,\overline{u}\,dx\ =\ 
k\int\limits_{Q_\infty}|u|^2dx\,. $$
Now we use that $u$ is also an eigenfunction and thus
$$ \int\limits_{Q_\infty}\bigl[|\nabla u|^2-k^2|u|^2\bigr]\,dx\ =\ 0\,. $$
Now we estimate (note that $|\sin\theta|<1$)
\begin{eqnarray*}
\Vert\nabla u\Vert^2_{L^2(Q_\infty)} & = & k^2\Vert u\Vert^2_{L^2(Q_\infty)}\ =\ 
k\biggl|\sin\theta\Im\int\limits_{Q_\infty}\frac{\partial u}{\partial x_1}\,\overline{u}\,dx
\biggr| \\ 
& < & k\Vert u\Vert_{L^2(Q_\infty)}\Vert\partial_1u\Vert_{L^2(Q_\infty)}\ =\ 
\Vert\nabla u\Vert_{L^2(Q_\infty)}\Vert\partial_1u\Vert_{L^2(Q_\infty)}\,,
\end{eqnarray*}
that is, $\Vert\nabla u\Vert_{L^2(Q_\infty)}<\Vert\partial_1u\Vert_{L^2(Q_\infty)}$ which 
implies that $\nabla u$ vanishes identically. Therefore $u$ is constant and thus zero by the 
boundary condition. \qed
\medlf
Now, applying Lemma \ref{sipe} we conclude that the unique solution $v(\epsilon)$ to 
\eqref{eq7} converges to $v$ in $X$ and the limiting function $v$ fulfils the equations
\ben
L_k v\ =\ f^{(k)}\quad\mbox{and}\quad PL^\prime(0)v\ =\ 0\,.
\enn 
The second equation provides an additional constraint on $v\in X$, that is (see \eqref{PK}), 
\be\label{or-v}
0\ =\ \int\limits_{Q_\infty}\bigl[\sin\theta\,\frac{\partial v}{\partial x_1}\,
\overline{\psi}\,dx\ -\ ik\cos^2\theta v\,\overline{\psi}\bigr]\,dx\quad\mbox{for all}\quad 
\psi\in \mathcal{N}.
\en
Setting $u=e^{ik\sin\theta x_1}v$ and $\phi=e^{ik\sin\theta x_1}\psi$, we return to 
quasi-periodic settings to get 
\ben
\sin\theta \int\limits_{Q_\infty}\frac{\partial u}{\partial x_1}\,\overline{\phi}\,dx\ =\ 
ik\int\limits_{Q_\infty}u\;\overline{\phi}\,dx\,
\enn
that is,
\be\label{orthogonal}
\int\limits_{Q_\infty}
\left(\sin\theta 
\frac{\partial u}{\partial x_1}-iku\right)\,\overline{\phi}\,dx\ =0\,\quad \mbox{for all}\quad\phi\in X_j\,,
\en

where $X_j\subset H^1_{\alpha,0}(Q_h)$ for all $h>R$ denotes the eigenspace \eqref{X_j} 
corresponding to the propagative wave number $\hat{\alpha}_j=k\sin\theta$. If we make the 
ansatz \eqref{unpert} for $u$: 
\begin{equation}\label{unpert}
u\ =\ u_0\ +\ \sum_{\ell=1}^{m_j}c_\ell\,\hat{\phi}_{\ell, j}\quad\mbox{in }D\,,\quad 
c_\ell\in \C\,,
\end{equation} 
where $u_0\in H^1_{\alpha,0}(Q_h)$ for all $h>h_0$ is a particular solution (for instance, given by Lemma \ref{lem3.4} (i)), then it follows from 
\eqref{orthogonal} that
\ben
\sum_{l=1}^{m_j}\left[\sin\theta\int\limits_{Q_\infty}\frac{\partial \hat{\phi}_{\ell, j}}{\partial x_1}\,\overline{\phi}\,dx-ik\int\limits_{Q_\infty}\hat{\phi}_{\ell, j}\;\overline{\phi}\,dx\right] c_l
=\int\limits_{Q_\infty}
\left(\sin\theta 
\frac{\partial u_0}{\partial x_1}-iku_0\right)\,\overline{\phi}\,dx
\enn
for all $\phi\in X_j$. Therefore, the 
coefficients $c_\ell$ should fulfill the finite-dimensional algebraic system
$$ (A-B)\,\mathcal{C}\ =\ Y $$
with 
\ben
& & \mathcal{C}\ =\ (c_1,c_2\cdots, c_{m_j})^\top\in\C^{m_j\times 1}\,,\quad 
A\ =\ \operatorname{diag}(a_{\ell,\ell})\in \C^{m_j\times m_j}\,, \\ 
& & B\ =\ (b_{\ell,\ell'})_{\ell,\ell'=1}^{m_j}\in \C^{m_j\times m_j},\quad 
Y\ =\ (y_1, y_2,\cdots, y_{m_j})^\top\in \C^{m_j\times 1}
\enn 
given by
\ben
y_\ell & := & \int\limits_{Q_\infty}
\left(\sin\theta 
\frac{\partial u_0}{\partial x_1}-iku_0\right)\,\overline{\phi_{l,j}}\,dx\,, \\
b_{\ell,\ell'} & := & i k\int\limits_{Q_\infty}\hat{\phi}_{\ell',j}\,
\overline{\hat{\phi}_{\ell,j}}\,dx\,, \\ 
a_{\ell,\ell}\ &:=&i/2\ \sin\theta\,\lambda_{\ell,j},
\enn
for $\ell,\ell'=1,2,\cdots, m_j$. Note that in deriving the entries of $A$ we have used the normalizations (see 
\eqref{evp})
\ben
\int\limits_{Q_\infty}\frac{\partial\hat{\phi}_{\ell,j}}{\partial x_1}\,
\overline{\hat{\phi}_{\ell',j}}\,dx\ =i\ \frac{\lambda_{\ell,j}}{2}\;\delta_{\ell,\ell'}\,.
\enn
\medlf
Well-posedness of scattering from unperturbed and perturbed periodic curves of Dirichlet kind 
is summarized as follows.
\begin{theorem}\label{TH-LAP} Let $k>0$ be fixed and $\theta\in(-\pi/2, \pi/2)$ be an arbitrary angle. Set 
$\alpha=k\sin\theta$ and suppose that $\alpha+n\neq\pm k$ for any $n\in \Z$ (that is, $\alpha$ is not a cut-off value)
\begin{itemize}
\item[(i)] In the unperturbed case, there exists a unique solution $u_{unpert}\in 
H^1_{\alpha,loc,0}(D)$ such that $u^{sc}_{unpert}:=u_{unpert}-u^{in}$ satisfies the upward 
$\alpha$-quasiperiodic  Rayleigh expansion \eqref{exc:b} as well as the constraint condition 
\be\label{OC}
\;\int\limits_{Q_\infty}\left(\sin\theta\frac{\partial u_{unpert}}{\partial x_1}-ik u_{unpert}\right)\,
\overline{\hat{\phi}_{\ell,j}}\,dx\ =0
\en
for all $\ell=1,2,\cdots, m_j$,
if $k\sin\theta=\hat{\alpha}_j$ is a propagative wave number. 
\item[(ii)] For locally perturbed periodic curves, there exists a unique solution 
$u\in H^1_{\alpha,loc,0}(\tilde{D})$ such that $u=u^{in}+u^{sc}_{unpert}+u^{sc}_{pert}$ 
in $\Sigma_R$ where $u^{sc}_{unpert}$ is given by assertion (i) and $u^{sc}_{pert}\in 
H_{loc}^1(\Sigma_R)$ fulfills the open waveguide radiation condition of Definition \ref{d-RC}
and the radiating part of $u^{sc}_{pert}$ satisfies the Sommerfeld radiation conditions of Definitions \ref{src} and \ref{src-i}.
\end{itemize}
\end{theorem}
\begin{proof}
(i) By Lemma \ref{lem3.4} (i), existence of $u_{unpert}\in H_{\alpha, loc,  0}(D)$ follows 
from the Fredholm alternative and uniqueness holds true if $\alpha$ is not a propagative wave number. 
If $\alpha=\hat{\alpha}_j$ for some $j\in J$,  we assume there are 
two solutions $u_{unpert}^{(1)}$ and $u_{unpert}^{(2)}$. Set 
$w=u_{unpert}^{(1)}-u_{unpert}^{(2)}$. It then follows from the limiting absorption argument  that the periodic function 
$v=e^{-i\alpha x_1}w\in \mathcal{N}$ fulfills the relation \eqref{or-v}, that is, 
$\langle PK'(0)v, \psi \rangle=0$ for all $\psi\in\mathcal{N}$.  Applying Lemma \ref{lem3.6} 
yields $v=0$ and thus $w=v e^{i\alpha x_1}=0$. 
\smalf
(ii) Once the unperturbed scattering problem is uniquely solvable, 
the uniqueness and existence of $u^{sc}_{pert}$ can be justified in the same way as 
\cite[Theorem 4.7]{HK22}.
\end{proof}


 As a corollary of Theorem  \ref{TH-LAP} (i), we obtain well-posedness of the following quasi-periodic boundary value problem:
\be\label{bvp-v}
\Delta v+k^2 v=0\quad\mbox{in}\quad D,
\qquad v=-g\quad\mbox{on}\quad \Gamma
\en
where $g=w|_\Gamma\in H^{1/2}_\alpha(\Gamma)$ with the function $w$ of the form
\ben
w(x)=\sum_{n\in \Z:\, |\alpha_n|<k} c_n\, e^{i(\alpha+n) x_1-i\beta_n x_2},\quad\quad c_n\in \C, \quad x\in Q_\infty.
\enn
\begin{corollary} Let $\alpha\in \R$ be arbitrary. The quasi-periodic boundary value problem \eqref{bvp-v} always admits a unique solution $v\in H^1_\alpha(Q_h)$ for all $h>h_0$ which fulfills the Rayleigh expansion condition $\eqref{exc:b}$. In the case that $\alpha=\hat{\alpha}_j$ is a propagative wavenumber, the unique solution $v$ is additionally  required to satisfy the orthogonal relation 
\be\label{OC}
\int\limits_{Q_\infty}\left(\sin\theta\frac{\partial (v+w)}{\partial x_1}-ik (v+w)\right)
\overline{\hat{\phi}_{\ell,j}}\,dx\ =0
\en
for all $\ell=1,2,\cdots, m_j$.

\end{corollary}

\begin{remark}
It remains unclear to us the well-posedness of the boundary value problem \eqref{bvp-v} with a general $g\in H_\alpha^{1/2}(\Gamma)$. For instance, $g$ is the restriction to $\Gamma$ of the incoming surface wave $e^{i \alpha_n x_1-i\beta_n x_2}$ with $|\alpha_n|>k$. In such a case, the function $f^{(k)}$ on the right hand of the variational formulation \eqref{eq7}, which can be expressed as
\ben
\langle f^{(k)},\psi\rangle := -2i\beta_n \,e^{-i\beta_n h}\int\limits_0^{2\pi}
\overline{\psi(x_1,h)}e^{in x_1}\,dx_1\quad\mbox{for all}\quad \psi\in X,
\enn
does not belong to the range of $L^{(k)}$. In fact, $f^{(k)}$ is not orthogonal to the null space $\mathcal{N}$ of $L^{(k)}$. However, if $\Gamma$ is given by a Lipschitz graph, it is well-known that the boundary value problem \eqref{bvp-v} admits a unique solution satisfying the $\alpha$-quasiperiodic upward Rayleigh expansion condition.
\end{remark}
\begin{remark}
For plane wave incidence, the approach of using the limiting absorption principle presented in this subsection also applies to the Neumann boundary condition  as well as transmission conditions for penetrable gratings. 
\end{remark}

\biglf

\subsection{Method of Approximation by Point Sources}\label{APS}
In this section we provide another proof of Theorem \ref{TH-LAP} by approximating a plane wave with point source waves. We shall prove that, when the location of the source tend to infinity, the total fields excited by point sources  converge to the total field of a plane wave and the limiting solution fulfills the same orthogonal constraint condition \eqref{OC}  at a propagative wavenumber. 

\smalf
We first consider the unperturbed scattering problem.

\begin{theorem}\label{TH-POI}
Let assumptions \ref{assump1} and \ref{assump2} hold and write $\hat{\theta}=(\sin\theta, -\cos\theta)^\top\in \R^2$ with a fixed $\theta\in(-\pi/2, \pi/2)$. Assume that 
\begin{itemize}
\item[(i)] $\alpha:=k\sin\theta$ is not a cut-off value in the sense of Def. \ref{d-exceptional} (i). 
\item[(ii)] The function $v(\cdot; \hat{\theta})\in H^1_{\alpha, loc, 0}(D)$ given by Theorem \ref{TH-LAP} (i) is the unique solution to the unperturbed scattering problem corresponding to the plane wave $u^{in}(x;\hat{\theta})=e^{ikx\cdot\hat{\theta}}$.
\end{itemize} 
Let $G_t=G(\cdot; z_t )$ with $z_t :=-t\hat{\theta}=(-t\sin\theta, t\cos\theta)$ be the unique total field of the unperturbed scattering problem of the point source at $z_t $ for $t \cos\theta>2h_0$ (see Proposition \ref{wps}). Then we have the convergence
\be\label{convergence1}
\frac{1}{\gamma} \lim_{t\rightarrow\infty} \left[\sqrt{t} e^{-ikt} G_t(x)\right]=v(x; \hat{\theta})\qquad\mbox{in}\quad H^1(Q_h),\quad \gamma:=\frac{e^{i\pi/4}}{\sqrt{8k\pi}},
\en
for any $h>h_0$.
\end{theorem}

\begin{remark}
The limiting function in \eqref{convergence1} relies essentially on the form of the unique solution $v(\cdot;\hat{\theta})$ to the unperturbed scattering problem, if $k\sin\hat{\theta}$ happens to be a propagative wavenumber. In this paper $v$ is derived from the LAP for approximating wave numbers. However, the analytical continuation arguments with respect to $\alpha$ or the LAP for approximating the refractive index in a slab gives a limiting solution satisfying different constraint conditions than \eqref{OC}; see \cite{K22p}. If $k\sin\hat{\theta}$ is not a propagative wavenumber, the limiting solutions obtained from these different approximation arguments are identical.
\end{remark}

\begin{proof} We carry out the proof following the lines in the proof of \cite[Theorem 5.2]{K22p} for inhomogeneous periodic layers. The proof will be divided into four steps.
\smalf

{\bf Step 1:} Reduction to the convergence proof for part of the radiating part. 

 As done in the proof of Theorem \ref{sy} (ii), for  each $z_t $ we choose a $t$-dependent cut-off function $\chi_t\in C_0^\infty(\real^2)$ with $\chi_t(x)=1$ 
for $|x-z_t |<\epsilon/2$ and $\chi_t(x)=0$ for $|x-z_t |\geq\epsilon$, where $\epsilon>0$ is fixed.
Then 
$w_t:=G_t-\chi_t\,\Phi(\cdot ; z_t )\in H^1_{loc}(D)$ coincides with $G_t$ for 
$|x-z_t |\geq\epsilon$ and satisfies $\Delta w_t+k^2w_t=-f_t$ in $D$ and $w_t=0$ on $\Gamma$ where
$$ f_t\ :=\ \Delta\chi_t\,\Phi(\cdot,z_t )+2\nabla\chi_t\cdot\nabla_x\Phi(\cdot,z_t )=(\Delta+k^2) [(\chi_t-1)\Phi(\cdot, z_t )]\in L^2(D) $$
has compact support. Let $w_{t, rad}$ and $w_{t, prop}$ be the radiating and propagating parts of $w_t$, respectively.
The radiating part $w_{t, rad}$ solves the inhomogeneous Helmholtz equation
\be\label{wr}
(\Delta +k^2)w_{t, rad}=-f_t-g_t\quad\mbox{in}\quad D, \qquad w_{t, rad}=0\quad\mbox{on}\quad \Gamma,
\en
where
\begin{eqngroup} 
\ben
g_t:=(\Delta +k^2)w_{t, prop}=\sum_{j\in J}\sum_{l=1}^{m_j} a_{l, j}(t)\,\varphi_{l,j},
\enn
\begin{equation} \label{eq-source1:b}
\varphi_{\ell,j}(x)\ =\ \left\{\begin{array}{cl} 2\,\psi_+^\prime(x_1)\,
\frac{\partial\hat{\phi}_{\ell,j}(x)}{\partial x_1}+\psi_+^{\prime\prime}(x_1)\,
\hat{\phi}_{\ell,j}(x) & \mbox{if }\lambda_{\ell,j}>0\,, \\ 
2\,\psi_-^\prime(x_1)\,\frac{\partial\hat{\phi}_{\ell,j}(x)}{\partial x_1}+
\psi_-^{\prime\prime}(x_1)\,\hat{\phi}_{\ell,j}(x) & \mbox{if }\lambda_{\ell,j}<0\,. 
\end{array}\right.
\end{equation}\end{eqngroup}
Note that $g_t$ is supported in the $x_1$-direction and exponentially decays in the $x_2$-direction and that the well-posdness of $w_{t, rad}$ is a consequence of \cite[Theorem 4.5]{HK22}. The coefficients 
$a_{l, j}(t)\in \C$ of the propagating part $w_{t, prop}$  have been computed explicitly in Theorem \ref{sy}  (ii), given by (see \eqref{alj})
\ben
a_{l, j}(t)=\frac{2\pi i}{|\lambda_{\ell,j}|}\,\overline{\hat{\phi}_{\ell,j}(z_t )}\quad \mbox{for all}\quad j\in J,\; l=1,2,\cdots, m_j.
\enn
This implies that
\be\label{es-a}
|a_{l, j}(t)|\leq c\, e^{-\delta t}, \quad ||g_t||_{L^\infty(D)}\leq c\, e^{-\delta t}\quad\mbox{for all}\quad t \geq 2h_0/\cos\theta,
\en with some $c>0$ independent of $t$. The same estimate holds for the propagating part (see \eqref{u2}):
\be\label{wp}
||w_{t, prop} ||_{L^\infty(D)}\leq C\,\sum_{j\in J}\sum_{l=1}^{m_j} |a_{l, j}(t)| \leq C\, e^{-\delta t}\quad\mbox{for all}\quad t \geq 2h_0/\cos\theta.
\en
The form of $w_t$ leads to a decomposition of $G_t$ as follows:
\ben
G_t=w_t+\chi_t\Phi(\cdot\,; z_t )=w_{t, rad}+w_{t, prop}+\chi_t\Phi(\cdot\, ; z_t ).
\enn Hence, the radiating part $G_{t, rad}$ of $G_t$ equals to $w_{t, rad}+\chi_t\Phi(\cdot\,; z_t )$, while the propagating part 
$G_{t, prop}$ coincides with $w_{t, prop}$. By $\eqref{wp}$ and the definition of $\chi_t$, 
\ben
\sqrt{t} e^{-ikt} || w_{t, prop}(\cdot)+\chi_t\Phi(\cdot\, ; z_t )||_{H^1(Q_h)}\rightarrow 0\quad\mbox{as}\quad t\rightarrow\infty,
\enn for any fixed $h>h_0$.
Therefore, it remains to consider the convergence of $w_{t, rad}$ as $t\rightarrow\infty$.

\smalf

{\bf Step 2:} Floquet-Bloch transform $w_{t, rad}$ to a family of quasi-periodic problems. 

For $g\in C_0^\infty(\R)$, the Floquet-Bloch transform $F$ is defined by
$$ (Fg)(x_1,\alpha)\ :=\ \sum_{n\in\Z} g(x_1+2\pi n)\,e^{-i2\pi n \alpha}\,,\qquad 
x_1\in\R,\ \alpha\in[-1/2, 1/2]\,. $$
The transform $F$ extends to an unitary operator from $L^2(\R)$ to 
$L^2((-1/2,1/2)\times(0, 2\pi))$. If $g$ depends on two variables $x_1$ and $x_2$ then the 
symbol $F$ means the Floquet-Bloch transform with respect to $x_1$. The inverse Floquet-Bloch transform is defined by $g=\int_{-1/2}^{1/2} (Fg)(\cdot, \alpha)d\alpha$. 
Taking the Floquet-Bloch transform on both sides of the equation \eqref{wr} yields
\be\label{eq:w}
(\Delta +k^2)w_{t, \alpha}=-(Ff_t)(\cdot, \alpha)-(Fg_t)(\cdot, \alpha)\quad\mbox{in}\quad Q_\infty, \qquad w_{t, \alpha}=0\quad\mbox{on}\quad \Gamma,
\en
where $w_{t, \alpha}=(F w_{t, rad})(\cdot, \alpha)\in L^2((-1/2, 1/2), H_{\alpha,0}^1(Q_\infty))$ (see \cite{L2017}). Here $H_{\alpha,0}^1(Q_\infty)$ is defined as the restriction of $H_{\alpha, loc, 0}^1(D)$ to $Q^\infty$.
 The above equation is understood in the variational sense that
\be\label{va-w}
\int\limits_{Q_\infty}[\nabla w_{t, \alpha}  \cdot \overline{\nabla \psi} - 
k^2 w_{t, \alpha}  \overline{\psi}]\,dx=\int\limits_{Q_\infty} [(Ff_t)(x, \alpha)+(Fg_t)(x, \alpha)]\overline{\psi}\,dx
\en
for all $\psi\in H^1_{\alpha,0}(Q_\infty)$.
We know from Theorem \ref{sy} (ii) and \cite[Theorem 3.5]{HK22} that for each $t>2h_0/\cos\theta$ this variational formulation is solvable for all $\alpha\in \R$ under the generalized Rayleigh expansion condition \eqref{angulaa-spectrum-rc} of $w_{t, \alpha}$, due to the orthogonality of the right hand of \eqref{eq:w} with the null space $X_{\alpha}$ by the choice of 
$a_{l,j}(t)$. Let $v_{t, \alpha}\in H^1_{\alpha}(U_{h_0})$ be the unique solution of the equation 
\ben
(\Delta +k^2)v_{t, \alpha}=-(Ff_t)(\cdot, \alpha)-(Fg_t)(\cdot, \alpha)\quad\mbox{in}\quad U_{h_0}, \qquad v_{t, \alpha}=0\quad\mbox{on}\quad \Gamma_{h_0},
\enn
together with the generalized Rayleigh expansion condition \eqref{angulaa-spectrum-rc} in $x_2>h_0$. It is easy to observe that 
\be\label{va-wv}
\int\limits_{U_{h_0}} [\nabla (w_{t, \alpha}-v_{t, \alpha})  \cdot \overline{\nabla \psi} - 
k^2 (w_{t, \alpha}-v_{t, \alpha})  \overline{\psi}]\,dx 
=\int\limits_{\Gamma_{h_0}} (\tilde{T}_k  w_{t, \alpha})\; \overline{\psi}\,ds
\en
for all $\psi\in H^1_{\alpha,0}(Q_\infty)$, where $\tilde{T}_k$ denotes the $\alpha$-quasiperiodic Dirichlet-to-Neumann map defined by \eqref{dtn2}. Simple calculations using \eqref{va-w} and \eqref{va-wv} show that the variational equation for $w_{t, \alpha}$ can be equivalently written as
\be \nonumber
\left(\mathbb{L}_\alpha  w_{t, \alpha}, \psi \right)_{H^1(Q_{h_0})}&:=&\int\limits_{Q_{h_0}}[\nabla w_{t, \alpha}  \cdot \overline{\nabla \psi} - 
k^2 w_{t, \alpha}  \overline{\psi}]\,dx-\int\limits_{\Gamma_{h_0}} \tilde{T}_k w_{t, \alpha} \overline{\psi}\,ds \\ \label{va-3} &=&\int\limits_{Q_{h_0}}(Fg_t)(x, \alpha)\overline{\psi}\,dx+\int\limits_{\Gamma_{h_0}} \frac{\partial v_{t, \alpha}}{\partial \nu} \overline{\psi}\,ds
\en
for all $\psi\in H^1_{\alpha,0}(Q_{h_0})$, which is defined as the restriction of $H^1_{\alpha,0}(Q_\infty)$ to $Q_{h_0}$. Note that by the choice of the cut-off function $\chi_t$ with $t\cos\theta>2h_0$, the function $f_t$ and thus $F f_t$ vanish in $Q_{h_0} $. Moreover,  we recall from \cite[Lemma 5.3]{K22p} that the normal derivative $\partial_\nu v_{t, \alpha}(x_1, h_0)$ can be computed explicitly as 
\ben
\frac{\partial v_{t, \alpha}(x_1, h_0)}{\partial x_2}&=&\frac{1}{2\pi}\sum_{l\in \Z} e^{i\sqrt{k^2-(l+\alpha)^2} (t\cos\theta-h_0)}\,e^{i(l+\alpha)(x_1+t\sin\theta)}\\
&&+\frac{1}{\sqrt{2\pi}}\sum_{l\in \Z}\int\limits_{h_0}^\infty (Fg_t)_l(y_2, \alpha) e^{i\sqrt{k^2-(l+\alpha)^2} (y_2-h_0)}\,dy_2\,e^{i(l+\alpha)x_1},
\enn
where $(Fg_t)_l(y_2, \alpha)$ are the Fourier coefficients of  $(Fg_t)_l(\cdot, y_2, \alpha)$, defined by
\ben
(Fg_t)_l(y_2, \alpha)=\frac{1}{\sqrt{2\pi}} \int\limits_0^{2\pi} (Fg_t)_l(y_1, y_2, \alpha) e^{-i (l+\alpha)y_1}dy_1.
\enn 
Hence, we get a family of quasi-periodic operator equations
\ben
\mathbb{L}_\alpha  w_{t, \alpha}=r_{t, \alpha}\quad \mbox{in}\quad H^1_{\alpha,0}(Q_{h_0}),
\enn
where $r_{t, \alpha}\in H^1_{\alpha,0}(Q_{h_0})$ is defined by
\ben
\left(r_{t, \alpha}, \psi \right)_{H^1(Q_{h_0})}&:=&\int\limits_{Q_{h_0}}(Fg_t)(x, \alpha)\overline{\psi}\,dx+\frac{1}{\sqrt{2\pi}}\sum_{l\in \Z} e^{i\sqrt{k^2-(l+\alpha)^2} (t\cos\theta-h_0)+i(l+\alpha)t\sin\theta}\;\overline{\psi_l(h_0)}\\
&&+\sum_{l\in \Z} \overline{\psi_l(h_0)} \int\limits_{h_0}^\infty (Fg_t)_l(y_2, \alpha) e^{i\sqrt{k^2-(l+\alpha)^2} (y_2-h_0)}\,dy_2.
\enn By the definition of $g_t$, $F$ and the estimate of $a_{l,j}(t)$ (see \eqref{es-a}), it follows that
\be\label{gt}
|(Fg_t)(x,\alpha)|+|\partial_\alpha (Fg_t)(x, \alpha)|\leq c\, e^{-\delta (t+|x_2|)}\quad\mbox{in}\quad x_2>h_0, 
\en for all $t>0$ and $\alpha\in[-1/2, 1/2]$. 

\smalf

{\bf Step 3:} Prove the convergence of the dominant part of $G_{t}$ as $t$ tends to infinity.

 Let $k=\tilde{l}+\kappa$ with $\hat{l}\in \Z$ and $\kappa\in(-1/2, 1/2]$. Then $\pm \kappa$ are cut-off values and they can decompose the interval $[-1/2, 1/2]$ into at most three open intervals $I_1\cup I_2\cup I_3$ such that their interiors are disjoint. Note that some of these intervals can degenerate into points and the cut-off values are contained in the boundary points of $I_m$, $m=1,2,3$.
Write  $k\sin\theta=\tilde{l}+\tilde{\alpha}$ with $\tilde{l}\in \Z$ and $\tilde{\alpha}\in(-1/2, 1/2]$.  Since $k\sin\theta$ is not a cut-off value, we suppose without loss of generality that $\tilde{\alpha}\in I_{\tilde{m}}$ is an interior point for some $\tilde{m}\in\{1,2,3\}$. Next, we find a subset $\mathcal{L}\subset \{-\tilde{l}, \cdots, \tilde{l}\}\subset \Z$ such that 
\ben
&&|\alpha+l|< k\qquad\mbox{for all}\quad \alpha\in I_{\tilde{m}},\quad  l\in \mathcal{L},\\
&&|\alpha+l|> k\qquad\mbox{for all}\quad \alpha\in I_{\tilde{m}},\quad  l\in \Z\backslash\mathcal{L}.
\enn

To find the dominant part of $G_{t}$, we decompose 
the Floquet-Bloch transform of the fundamental solution $\Phi_t(x)=\Phi(x;  z_t )$ with $t\cos\theta>2h_0$ and $x\in Q_{h_0}$ into 
\ben
\Phi_{t, \alpha}(x)=(F\Phi_t)(x, \alpha)&=&\frac{i}{4\pi}\sum_{l\in \Z} \frac{e^{i (\alpha+l) (x_1+t\sin\theta)+i\sqrt{k^2-(l+\alpha)^2} (t\cos\theta-x_2)}}{\sqrt{k^2-(l+\alpha)^2}}
=\Phi^{(1)}_{t, \alpha}+\Phi^{(2)}_{t, \alpha},
\enn
where $\Phi^{(2)}_{t, \alpha}:=\Phi_{t, \alpha}-\Phi^{(1)}_{t, \alpha}$ with
\ben
\Phi^{(1)}_{t, \alpha}(x):=\frac{i}{4\pi}\sum_{l\in \mathcal{L}} \frac{e^{i t [ (\alpha+l) \sin\theta+\sqrt{k^2-(l+\alpha)^2} \cos\theta]}}
{\sqrt{k^2-(l+\alpha)^2}} \; v^{in}_{l, \alpha}(x), \;\;v^{in}_{l, \alpha}(x):=e^{i \left( (\alpha+l) x_1-\sqrt{k^2-(l+\alpha)^2} x_2\right)}.
\enn
Note that the  Floquet-Bloch transform of $\Phi(x; y)$ is nothing else but the quasi-periodic fundamental solution to the Helmholtz equation.
For $l\in\mathcal{L}$, $v^{in}_{l, \alpha}$ is an incident plane wave with the unit direction  $(\alpha+l, -\sqrt{k^2-(l+\alpha)^2})/k$. We denote by $v_{l, \alpha}$ the unique $\alpha$-quasiperiodic  total field generated by $v^{in}_{l, \alpha}$ (see Theorem \ref{TH-LAP} (i)). In particular, $v_{\tilde{l}, \tilde{\alpha}}=v(\cdot;\hat{\theta})$ when $l=\tilde{l}$ and $\alpha=\tilde{\alpha}$.
We remark that $v_{l, \alpha}$ is required to satisfy the orthogonal condition  \eqref{OC}, if $\alpha$ is a propagative wavenumber. By linear superposition, the total field  excited by  $\Phi^{(1)}_{t, \alpha}$, which we denote by $w^{(1)}_{t, \alpha}$,  can be represented as 
\ben
w^{(1)}_{t, \alpha}(x)&=&\frac{i}{4\pi}\sum_{l\in \mathcal{L}\backslash\{\tilde{l}\}} \frac{e^{i t [ (\alpha+l) \sin\theta+\sqrt{k^2-(l+\alpha)^2} \cos\theta]}}
{\sqrt{k^2-(l+\alpha)^2}} \; v_{l, \alpha}(x)+ W_{t, \alpha}(x),\\
W_{t, \alpha}(x)&:=&\frac{i}{4\pi}\frac{e^{i t [ (\alpha+\tilde{l}) \sin\theta+\sqrt{k^2-(\tilde{l}+\alpha)^2} \cos\theta]}}
{\sqrt{k^2-(\tilde{l}+\alpha)^2}} \; v_{\tilde{l}, \alpha}(x).
\enn
It was proved in \cite{K22p} that the inverse Floquet-Bloch transform of  $w^{(1)}_{t, \alpha}$ (more precisely, $W_{t, \alpha}$)  constitutes  the dominant part of $G_{t}$ as $t\rightarrow\infty$.
In fact,  using stationary arguments one deduces that (see e.g. \cite[Section 5]{K22p}) 
\ben
\int\limits_{I_{\tilde{m}}}W_{t, \alpha}(x)\,d\alpha=\gamma \frac{e^{itk}}{\sqrt{t}} v(x; \hat{\theta})+o(t^{-1/2}),
\enn
and using partial integration yields (see Subsection \ref{asy} in the Appendix)
\be\label{partial}
\int\limits_{I_{\tilde{m}}} [w^{(1)}_{t, \alpha}(x)-W_{t, \alpha}(x)]\,d\alpha=O(t^{-1}),\quad
\int\limits_{[-1/2, 1/2]\backslash I_{\tilde{m}}}w^{(1)}_{t, \alpha}(x)\,d\alpha=O(t^{-1}),
\en
as $t\rightarrow\infty$. This proves
\ben
 \lim_{t\rightarrow\infty} \left[\sqrt{t} e^{-ikt} \int\limits_{-1/2}^{1/2} w^{(1)}_{t, \alpha}  (x)\,d\alpha  \right]=v(x; \hat{\theta})\qquad\mbox{in}\quad H^1(Q_{h}),\;h>h_0.
\enn
\smalf
{\bf Step 4:} Show the decay of the remaining part. 

To prove \eqref{convergence1}, we only need to show for $w^{(2)}_{t, \alpha}:= w_{t, \alpha}-w^{(1)}_{t, \alpha}$ that
\be\label{limit}
\lim_{t\rightarrow\infty}\left[ \sqrt{t} e^{-ikt} \int\limits_{-1/2}^{1/2}w^{(2)}_{t, \alpha}  (x) \,d\alpha\right]=0\quad\mbox{in}\quad H^1(Q_{h}),\; h>h_0.
\en
Recalling the variational formulation for the total field $w^{(1)}_{t, \alpha}$ (see \eqref{va-3} for the definition of $\mathbb{L}_\alpha$), 
\ben
\left(\mathbb{L}_\alpha w^{(1)}_{t, \alpha}, \psi \right)_{H^1(Q_{h_0})}:=\frac{1}{\sqrt{2\pi}}\sum_{l\in \mathcal{L}} e^{i\sqrt{k^2-(l+\alpha)^2} (t\cos\theta-h_0)+i(l+\alpha)t\sin\theta}\;\overline{\psi_l(h_0)}
\enn for all $\psi\in H^1_{\alpha,0}(Q_{h_0})$,
we find that $w^{(2)}_{t, \alpha}$ are solutions of $\mathbb{L}_\alpha  w^{(2)}_{t, \alpha}=\sum_{j=1}^3 r^{(j)}_{t, \alpha}$, where 
\ben
&&\left(r^{(1)}_{t, \alpha}, \psi \right)_{H^1(Q_{h_0})}:= \int\limits_{Q_{h_0}}(Fg_t)(x, \alpha)\overline{\psi}\,dx,\\
&&\left(r^{(2)}_{t, \alpha}, \psi \right)_{H^1(Q_{h_0})}:=\frac{1}{\sqrt{2\pi}}\sum_{l\notin \mathcal{L}} e^{-\sqrt{(l+\alpha)^2-k^2} (t\cos\theta-h_0)+i(l+\alpha)t\sin\theta}\;\overline{\psi_l(h_0)},\\
&&\left(r^{(3)}_{t, \alpha}, \psi \right)_{H^1(Q_{h_0})}:= \sum_{l\in \Z} \overline{\psi_l(h_0)} \int\limits_{h_0}^\infty (Fg_t)_l(y_2, \alpha) e^{i\sqrt{k^2-(l+\alpha)^2} (y_2-h_0)}\,dy_2.
\enn
Since every cut-off value is assumed to be no propagative wavenumber, one may divide the interval [-1/2, 1/2] into the union $\Lambda_1\cup \Lambda_2$ of two types of closed sub-intervals with non-intersecting interiors, where $\Lambda_1$ does not contain any propagative wavenumber and $\Lambda_1$ contains no cut-values. In $\Lambda_1$, one can deduce from the decaying of $Fg_t$ (see \eqref{gt})  and partial integration that (see \cite{K22p} for the details)
\ben
\int\limits_{\Lambda_1} ||w^{(2)}_{t, \alpha} ||_{H^1(Q_{h_0})}\,d\alpha\leq c\, \left(\sum_{l\notin \mathcal{L}} 
\int\limits_{\Lambda_1} e^{-t\sqrt{(l+\alpha)^2-k^2}\cos\theta}d\alpha\right)^{1/2}\leq c\,t^{-1}
\enn Since $r_{t, \alpha}^{(j)}$ are differentiable with respect to $\alpha\in \Lambda_2$ for $j=1,2,3$, the integral over $\Lambda_2$ can be estimated by applying Lemma \ref{sipe} to get (see also \cite{K22p})
\ben
\int\limits_{\Lambda_2} ||w^{(2)}_{t, \alpha} ||_{H^1(Q_{h_0})}\,d\alpha\leq c\,t\, e^{-\delta t}\qquad\mbox{for all}\quad t\cos\theta\geq 2h_0.
\enn
Combining the previous two estimates yields \eqref{limit} and thus finishes the proof of  Theorem \ref{TH-POI}.
\end{proof}

\smalf

Now we study the limit of the Green's function to the locally perturbed scattering problem when the source position tend to infinity. 
\begin{theorem}\label{TH-PL}
Let Assumptions \ref{assump1}, \ref{assump2} and \ref{assump3} hold and write $\hat{\theta}=(\sin\theta, -\cos\theta)^\top\in \R^2$ with $\theta\in(-\pi/2, \pi/2)$. Assume that $\alpha:=k\sin\theta$ is not a cut-off value in the sense of Def. \ref{d-exceptional} (i). Let $u_t=u(\cdot; z_t )$ with $z_t :=-t\hat{\theta}$ be the unique total field of the perturbed scattering problem of the point source at $z_t $ for $t \cos\theta>2h_0$ (see Proposition \ref{wps}), which satisfies the open waveguide radiation condition of Definition \ref{d-RC}.
Then we have the convergence
\be\label{convergence}
\frac{1}{\gamma} \lim_{t\rightarrow\infty} \left[\sqrt{t} e^{-ikt} u_t(x)\right]=w(x; \hat{\theta})
\qquad\mbox{in}\quad H^1(\tilde{D}_R),\quad \gamma:=\frac{e^{i\pi/4}}{\sqrt{8k\pi}},
\en
for any $R>\pi$, where $w\in H^1_{loc, 0}(\tilde{D})$ with the decomposition $w=u^{in}+u^{sc}_{unpert}+u^{sc}_{pert}$ in $\Sigma_R$ denotes the unique solution to the perturbed scattering problem corresponding to the plane wave $u^{in}(x;\hat{\theta})=e^{ikx\cdot\hat{\theta}}$ specified in Theorem \ref{TH-LAP} (ii). 
\end{theorem}
\begin{remark}
In the absence of the defect, $u_t$ coincides with the Green's function $G_t$ to the scattering problem in perfectly periodic structures and $w=u^{in}+u^{sc}_{pert}$ coincides with the limiting function $v$ specified in  Theorem \ref{TH-POI}.
\end{remark}
\begin{proof}
By the proof of Proposition \ref{wps} (see \cite{HK22}), the total field $u_t$ can be decomposed into $u_t=G_t+ u^{sc}_{t, pert}$ in $\Sigma_R$ where $G_t$ is the Green's function to the unperturbed scattering problem and $u^{sc}_{t, pert}$ corresponding to the defect satisfies the open waveguide radiation condition. It follows from Theorem \ref{TH-POI} that
\be\label{con}
\frac{1}{\gamma} \lim_{t\rightarrow\infty} \left[\sqrt{t} e^{-ikt}\; G_t\right]=u^{in}+u^{sc}_{unpert}=:w^{in}
\qquad\mbox{in}\quad H^1(Q_h)
\en
for all $h>h_0$. To prove the convergence \eqref{convergence}, we define
\ben
v_t:=1/\gamma \sqrt{t} e^{-ikt} u_{t}-w\quad\mbox{in}\quad \tilde{D},
\enn
which can be considered as the total field corresponding to $v^{in}_t:=1/\gamma \sqrt{t} e^{-ikt} G_{t}-w^{in}$. It is obvious that $v_t-v_t^{in}$ fulfills the open waveguide radiation condition. 
\smalf

Choose $R>\pi$  such that there is no bound state to the Helmholtz equation over the domain $\Sigma_R$ and that $k^2$ is not the Dirichlet eigenvalue of the nagative Laplacian operator over $\tilde{D}_R$. We suppose without loss of generality that the domain $\tilde{D}_R$ is Lipschitz. Otherwise, one can slightly change the shape of $C_R$ to get a Lipschitz domain.
On  the artificial curve $C_R$, one may construct the Dirichlet-to-Neumann operator $\Lambda$ that is equivalent to the open waveguide radiation condition. The operator $\Lambda: H^{1/2}_0(C_R)\rightarrow H^{-1/2}(C_R)$ has been proved to be bounded  and $-\Lambda$ can be decomposed into the sum of coercive operator and a compact operator; see \cite[Lemma 3.9]{HK22}. With the aid of this DtN operator, one deduces the following boundary value problem for $v_t\in\{u\in H^1(\tilde{D}_R): u=0\;\mbox{on}\; \partial \tilde{D}_R\cap \tilde{\Gamma}\}$:
\ben
({\rm BVP}):\quad\left\{\begin{array}{lll}
(\Delta +k^2) v_t=0\quad&&\mbox{in}\quad \tilde{D}_R, \\
\partial_\nu v_t=\Lambda\, v_t +(\partial_\nu v^{in}_t -\Lambda v^{in}_t) \quad&&\mbox{on}\quad C_R,
\end{array}\right.
\enn
where $\nu$ denotes the normal direction at $C_R$ pointing into $\Sigma_R$. The well posedness of the above boundary value problem follows from mapping properties of the DtN operator together with the assumption that there is no bound state over $\tilde{D}$ (see \cite[Theorem 2.9 (ii)]{HK22}). Hence, using \eqref{con} and  the boundedness of $\Lambda$ we arrive at
\ben
||v_t||_{ H^1(\tilde{D}_R)}\leq c\, || \partial_\nu v^{in}_t -\Lambda v^{in}_t ||_{H^{-1/2}(C_R)}\leq c\, || v^{in}_t  ||_{H^{1}(D_R)}
\rightarrow 0,
\enn
as $t \rightarrow \infty$, which proves \eqref{con}.

\end{proof}

\section{Uniqueness results to inverse scattering}\label{sec4}
This section is concerned with uniqueness in determining the shape and location of the defect $\tilde{\Gamma}\backslash \Gamma$ from near/far-field data incited by plane or point source waves at a fixed wavenumber. We suppose that the unperturbed grating profile $\Gamma=\partial D$ is {\em a priori} known with the period $2\pi$. 
Although we only discuss a localized defect appearing on the scattering interface, 
 the uniqueness results of this section carry over to a perturbation caused by a bounded Dirichlet obstacle embedded inside $D$. 

\subsection{Uniqueness With Infinitely Many Point Source Waves}\label{sub5-1}
Let $G(x; y)$ ($x\neq y$) be the total field (Green's function) to the perturbed scattering problem with $u^{in}=\Phi(x; y)$; see Proposition \ref{wps}.

\begin{theorem}\label{uni-sour-inf} Let $\tilde{\Gamma}$ be a local perturbation of the periodic curve $\Gamma$ and suppose that 
$\max\{x_2: x\in \Gamma\cup \tilde{\Gamma}\}<h$ for some $h\in \R$.  
Then $\tilde{\Gamma}$ can be uniquely determined by the near-field measurement data $\{G(x_1, h; y_j): x_1\in(a,b), y_j\in U_h, j=1,2,\cdots\}$, incited by infinitely many point source waves. 
\end{theorem}
\begin{proof}
Suppose that there are two local perturbations $\tilde{\Gamma}_1$ and $\tilde{\Gamma}_2$ which both lie below the line $x_2=h$. Denote by $G_{\ell}(x; y_j)$ ($\ell=1,2$) the total fields corresponding to $\tilde{\Gamma}_{\ell}$ and the incoming source wave $\Phi(x; y_j)$, and 
 let $\tilde{D}_j$ be the domain above $\tilde{\Gamma}_j$. 
 Assuming 
\be\label{h}
G_1(x_1, h; y_j)=G_2(x_1, h; y_j)\quad \mbox{for all}\quad x_1\in(a, b),\; j\in \N,
\en
we need to prove $\tilde{\Gamma}_1=\tilde{\Gamma}_2$.
By the analyticity of $G_\ell (\cdot ; y_j)$ on $x_2=h$, we deduce from \eqref{h} that $G_1(x; y_j)=G_2(x; y_j)$ on $x_2=h$ for all $j\in \N$.  With the open waveguide radiation condition of Def. \ref{d-RC}, there exists a unique solution to 
 the Dirichlet boundary value problem of the Helmholtz equation in the upper half-plane $x_2>h$; we refer to Lemma \ref{lem4.2} and Remark \ref{rem4.3} below for the proof. 
Hence, for fixed $j\in \N$, the functions $G_1(\cdot; y_j)$ and $G_2(\cdot; y_j)$ must coincide in $x_2>h$ and by unique continuation also coincide in $\Omega\backslash\{y_j\}$ where $\Omega$ denotes the unbounded component of $\tilde{D}_1\cap \tilde{D}_2$. Consequently, the total fields $G_{\ell}$ ($\ell=1,2$) vanish on $\partial \Omega$.
\medlf
 If $\tilde{\Gamma}_1\neq\tilde{\Gamma}_2$, we shall derive a contradiction as follows. Switching the notation if necessary,  we can  assume that (see Figure \ref{fig1})  
\ben
D^*=[ (\tilde{D}_1\cup \tilde{D}_2)\,\backslash\, \Omega\, ]\,\cap\, \tilde{D}_1\neq\emptyset. 
\enn
\begin{figure}[htb]
  \centering
    \includegraphics[width=0.75\textwidth]{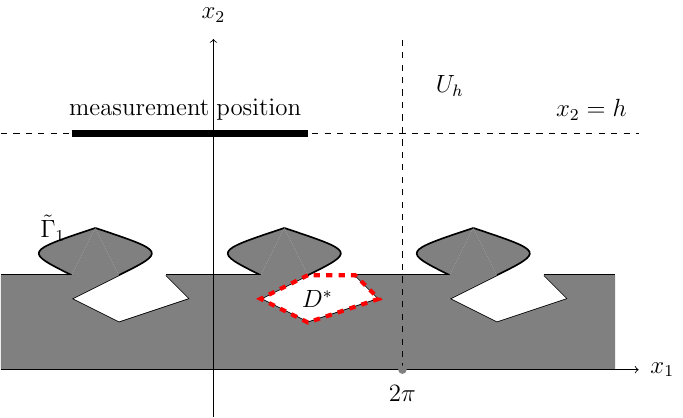}
\caption{Illustration of the gap domain $D^*\subset \tilde{D}_1\backslash\overline{\tilde{D}}_2$ between two local perturbations. Here $\tilde{\Gamma}_1=\Gamma$ is identical with the unperturbed grating curve.}\label{fig1}
\end{figure}

It is obvious that $\partial D^*\subset \tilde{\Gamma}_1 \cup (\tilde{\Gamma}_2\cap \partial \Omega)$. Noting that $y_j\notin D^*$ and $G_1=G_2=0$ on  $\tilde{\Gamma}_2\cap \partial \Omega$, we obtain
\ben
(\Delta +k^2)G_1(x; y_j)=0\quad\mbox{in}\quad D^*,\quad G_1(\cdot; y_j)=0\quad\mbox{on}\quad \partial D^*
\enn for all $j\in \N$. This implies that there exist infinitely many Dirichlet eigenfunctions $G_1(\,\cdot\,; y_j)$ for the negative Laplacian operator over the bounded domain $D^*$ with the eigenvalue $k^2$. Now, it suffices to show the linear independence of 
$G_1(\,\cdot \,; y_j)$, which together with the finite dimensional Dirichlet eigenspace (irrespective of boundary regularities) could lead to a contradiction. Assume that
\ben
\sum_{j=1}^M \lambda_j\, G_1(x; y_j)=0,\qquad x\in D^*,
\enn
for some constants $\lambda_j\in \C$, where $y_j\in U_h$ for $j=1,2,\cdots, M$ are distinct point sources.  Since $D^*\subset \tilde{D}_1$, applying the unique continuation yields
\be\label{G1}
\sum_{j=1}^J \lambda_j\, G_1(x; y_j)=0,\qquad\mbox{for all}\quad x\in \tilde{D}_1\backslash\{y_j\}_{j=1}^M. 
\en
Now, letting $x\rightarrow y_j$ in \eqref{G1} and using the boundedness of $G_1(y_j; y_\ell)$ for $\ell\neq j$, we obtain $\lambda_j=0$. The arbitrariness of $j=1,2,\cdots, J$ implies that the total fields corresponding different point sources are indeed linearly independent. This finishes the proof of $\tilde{\G}_1=\tilde{\G}_2$.
\end{proof}

In the proof of Theorem \ref{uni-sour-inf}, we need the following uniqueness result to  the homogeneous Dirichlet boundary value problem of the Helmholtz equation in the half plane $x_2>h$. Let $\sigma>0$ be the parameter of the cut-off function given in Def. \ref{d-RC}.
\begin{lemma}\label{lem4.2}
Let $u\in H^1_{loc}(U_h)$ be a solution of the Helmholtz equation $\Delta u+k^2u=0$ in $U_h$
such that $u=0$ on $x_2=h$. Furthermore, let $u$ be of the form $u=u_{rad}+u_{prop}$ where 
$u_{rad}\in H^1(U_h\setminus U_H)$ for all $H>h$ satisfies the generalized angular spectrum radiation 
condition \eqref{angulaa-spectrum-rc} and $u_{prop}=\sum_{j\in J}u_j$ where $u_j=\sum_{\ell=1}^{m_j}a_{\ell,j}^\pm
\hat{\phi}_{\ell,j}$ for $\pm x_1\geq\sigma$. Then $u$ vanishes in $U_h$.
\end{lemma}
\textbf{Proof:} Let $I\subset\real$ be any bounded interval. Set $I_n=\{t+2\pi n:t\in I\}$ for 
$n\in\ganz$. Then, for sufficiently large $n>0$,
\begin{eqnarray*}
\sum_{j\in J}e^{2\pi n\alpha_j i}\int\limits_I u_j(x_1,h)\,dx_1 & = &
\sum_{j\in J}\int\limits_I u_j(x_1+2\pi n,h)\,dx_1\ =\ 
\int\limits_{I_n} u_{prop}(x_1,h)\,dx_1 \\
& = & -\int\limits_{I_n} u_{rad}(x_1,h)\,dx_1\ \longrightarrow\ 0
\end{eqnarray*}
as $n\to\infty$. Here we have used the quasiperiodicity of $\hat{\phi}_{l,j}$ and the definition of $u_{prop}$.
Set $b_j:=\int_I u_j(x_1,h)\,dx_1$ for abbreviation. Then $\sum_{j\in J}b_j\,e^{2\pi n\alpha_j i}$
tends to zero as $n\to\infty$. By induction with respect to $|J|$ (number of elements) one 
proves that all $b_j$ vanish. Indeed, this is obviously true for $|J|=1$. Let it hold for 
$|J|=p$ and let $\hat{J}=J\cup\{\alpha_\ell\}$ with $|J|=p$ and $\alpha_\ell\notin
\{\alpha_j:j\in J\}$ and
\begin{equation} \label{aux}
\sum_{j\in J}b_j\,e^{2\pi n\alpha_j i}\ +\ b_\ell\,e^{2\pi n\alpha_\ell i}\ 
\longrightarrow\ 0\,,\quad n\to\infty\,.
\end{equation}
Multiplication of this formula by $e^{2\pi \alpha_\ell i}$ yields the first of the following 
formula:
\begin{eqnarray*}
\sum_{j\in J}b_j\,e^{2\pi n\alpha_j i}\,e^{2\pi\alpha_\ell i}\ +\ 
b_\ell\,e^{2\pi(n+1)\alpha_\ell i} & \longrightarrow & 0\,,\quad n\to\infty\,, \\
\sum_{j\in J}b_j\,e^{2\pi(n+1)\alpha_j i}\ +\ b_\ell\,e^{2\pi(n+1)\alpha_\ell i} & 
\longrightarrow & 0\,,\quad n\to\infty\,.
\end{eqnarray*}
Note that the second one is (\ref{aux}) for $n+1$ instead of $n$.
Subtraction of the previous two relations yields
$$ \sum_{j\in J}b_j\,e^{2\pi n\alpha_j i}\,[e^{2\pi\alpha_\ell i}-
e^{2\pi\alpha_j i}]\ \longrightarrow\ 0\,,\quad n\to\infty\,. $$
Now we apply the assumption of induction to $\tilde{b}_j=[e^{2\pi\alpha_\ell i}-
e^{2\pi\alpha_j i}]b_j$ which gives $b_j=0$ for all $j\in J$ and thus also $b_\ell=0$.
\smalf
Therefore, $\int_I u_j(x_1,h)\,dx_1=0$ for all $j\in J$ and all intervals $I$. This proves 
that $u_{prop}$ vanishes for $x_1>\sigma$. The same argument for $n\to-\infty$ yields that 
$u_{prop}$ vanishes for $x_1<-\sigma $. Therefore $u$ itself satisfies the generalized angular 
spectrum radiation condition and vanishes for $x_2=h$. This yields $u=0$ by arguing the same as in the proof of the last assertion of \cite[Appendix, Lemma 7.1]{K22}.
 \qed
\begin{remark}\label{rem4.3}
If $u$ vanishes on a locally perturbed periodic curve $\tilde{\Gamma}$ instead of the straight line $x_2=h$, it follows from \cite[Theorem 2.8]{HK22} that we still have
$u_{prop}=0$. However, $u=u_{rad}\in H^1_0(\tilde{D})$ becomes a bound state over the domain $\tilde{D}$;  see \cite[Theorem 2.9]{HK22}. The above lemma presents a simple proof for the vanishing of the propagating part when $\tilde{\G}=\partial \tilde{D}$ is a straight line. 
\end{remark}
\medlf

The proof of Theorem \ref{uni-sour-inf} does not carry over to the Neumann  boundary condition, because the property of a finite dimensional eigenspace in the Neumann case requires  boundary smoothness assumptions which usually cannot be fulfilled. 
Below we present another proof relying on the blowing up argument of \cite{Isa90, KK93}, which applies to  the Neumann and Impedance boundary conditions provided the well-posedness of forward scattering problems can be justified.

\medlf
\begin{theorem}\label{Th-inverse-source}
Under the assumption of Theorem  \ref{uni-sour-inf}, the locally perturbed defect $\tilde{\Gamma}$ can be uniquely determined by the near-field measurement data $\{G(x_1, h; y_1, h): x_1\in(a,b), y_1\in (c, d)\}$. Here $(a, b), (c, d)\subset \R$ are finite intervals without intersections.
\end{theorem}
\begin{proof}
We keep the notations in the proof of Theorem  \ref{uni-sour-inf} to obtain
\be\label{G12}
G_1(x; y_1, h)=G_2(x; y_1, h)\quad \mbox{for all}\quad x\in \Omega, \; y_1\in(c, d).
\en
Using the symmetry of $G_\ell(x; y)$ (see Theorem \ref{sy}), we deduce from \eqref{G12} and the unique continuation that 
\be\label{R3}
G_1^{sc}(x; y)=G_2^{sc}(x; y)\quad \mbox{for all}\quad x, y\in \Omega, \; x\neq y.
\en
If $\tilde{\Gamma}_1\neq\tilde{\Gamma}_2$, without loss of generality we can choose a point $y^*$ and a sub-boundary $\mathcal{S}$ of $\tilde{\Gamma}_1$ such that
$y^*\in \mathcal{S}\subset (\partial \Omega\cap \tilde{\Gamma}_1)\cap \tilde{D}_2$  and
\ben
y^{(m)}:=y^*+\nu(y^*)/m \in \Omega\cap \tilde{D}_2
\enn
for all $m\geq M$ with some $M\in \N$,
where $\nu(y^*)\in \mathbb{S}:=\{x\in \R^2: |x|=1\}$ denotes the unit normal direction at $y^*\in \tilde{\Gamma}_1$ pointing into $\tilde{D}_1$. Since $y^*$ is bounded away from $\tilde{\Gamma}_2$, well posdness of the forward scattering problem for $\tilde{\Gamma}_2$ implies that
\be\label{R1}
\lim_{m\rightarrow\infty}||G_2^{sc}(x; y^{(m)})||_{H^{1/2}(\mathcal{S})}=||G_2^{sc}(x; y^*)||_{H^{1/2}(\mathcal{S})}<\infty.
\en 
On the other hand, it follows from the Dirichlet boundary condition $G_1^{sc}(x; y^*)=-\Phi(x; y^*)$ on $\mathcal{S}$ that 
\be\label{R2}
\lim_{m\rightarrow\infty}||G_1^{sc}(x; y^{(m)})||_{H^{1/2}(\mathcal{S})}=||\Phi(x; y^*)||_{H^{1/2}(\mathcal{S})}=\infty,
\en
due the the singular behaviour $\Phi(x;y^*)=O(\ln |x-y^*|)$ as $|x-y^*|\rightarrow 0$.
The previous two relations \eqref{R1} and \eqref{R2} obviously contradict with the identity  \eqref{R3}. This contradiction proves that $\tilde{\Gamma}_1=\tilde{\Gamma}_2$.
\end{proof}

\subsection{Uniqueness With Infinitely Many Plane Waves} Let $u^{in}$ be a plane wave with fixed wavenumber $k>0$.
To specify the dependence on the incident angle $\theta\in (-\pi/2, \pi/2)$, we rewrite the unique total field $u\in H^1_{loc, \alpha, 0}(\tilde{D})$ to the perturbed scattering problem as (e.g., Theorem \ref{TH-LAP} (ii))
\be\label{total}
u(x;\theta)=u_{unpert}(x;\theta)+u^{sc}_{pert}(x;\theta)\qquad\mbox{in}\quad 
\Sigma_R 
\en
where $u_{unpert}(x;\theta)=u^{in}(x;\theta)+u^{sc}_{unpert}(x;\theta)\in H^1_{\alpha, loc, 0}(D)$ is the total field to the unperturbed scattering problem, and $u^{sc}_{pert}(x;\theta)\in 
H_{loc}^1(\Sigma_R)$ is caused by the local defect which fulfils the open waveguide radiation condition of Definition \ref{d-RC}.
Note that, if $k\sin\theta=\hat{\alpha}_j+n$ for some $n\in \Z$ and some propagative wavenumber $\hat{\alpha}_j$ $(j\in J)$, the unperturbed total field $u_{unpert}$ is supposed to fulfil the additional constraint of Theorems \ref{TH-LAP} and \ref{TH-POI}.

\begin{theorem}\label{uni-plane-inf} Let $\tilde{\Gamma}$ be a local perturbation of the periodic curve $\Gamma$ and suppose that 
$\max\{x_2: x\in \Gamma\cup \tilde{\Gamma}\}<h$ for some $h\in \R$.  
Then $\tilde{\Gamma}$ can be uniquely determined by the near-field measurement data $\{u(x_1, h; \theta_m): x_1\in(a,b),  m=1,2,\cdots\}$, incited by infinitely many plane waves with distinct incident angles $\theta_m\in (-\pi/2, \pi/2)$. 
\end{theorem}
\begin{proof} We shall carry out the proof by arguing analogously to the proof of Theorem \ref{uni-sour-inf}. It suffices to prove the linear independence of the total fields caused by different directions.

Set $\alpha(n):=k\sin\theta_n$ for $n=1,\ldots,N$. We recall that the total field $u_n$ 
corresponding to the incident angle $\theta_n$ has the decomposition into $u_n=u_n^{in}+
u^{sc}_{n,unpert}+u^{sc}_{n,pert}$ where $u_n^{in}(x)=e^{i\alpha(n)x_1-i\sqrt{k^2-\alpha(n)^2}x_2}$ and 
$u^{sc}_{n,pert}$ satisfies the open waveguide radiation condition and 
$u^{sc}_{n,unpert}$ is $\alpha(n)$-quasi-periodic; that is, has a Rayleigh expansion in the  form
$$ u^{sc}_{n,unpert}(x)\ =\ \sum_{\ell\in\ganz}u_{\ell,n}\,e^{i\sqrt{k^2-(\ell+\alpha(n))^2}x_2}\,
e^{i(\ell+\alpha(n))x_1}\,,\quad x_2\geq h\,. $$
Let now $\sum_{n=1}^N\lambda_nu_n=0$ in $\tilde{D}$. For fixed $m\in\{1,\ldots,N\}$, $R>0$ and 
$x_2>h$ we have
\begin{eqnarray} 
0 & = & \frac{1}{2R}\sum_{n=1}^N\lambda_n\int\limits_{-R}^Ru_n(x)\,e^{-i\alpha(m)x_1}\,dx_1  
\nonumber \\
& = & \sum_{n=1}^N\lambda_n\,e^{-i\sqrt{k^2-\alpha(n)^2}x_2}\frac{1}{2R}\int\limits_{-R}^R
e^{i(\alpha(n)-\alpha(m))x_1}\,dx_1 \nonumber \\
& & +\ \sum_{n=1}^N\sum_{\ell\in\ganz}\lambda_n\,u_{\ell,n}\,e^{i\sqrt{k^2-(\ell+\alpha(n))^2}x_2}\,
\frac{1}{2R}\int\limits_{-R}^R e^{i(\ell+\alpha(n)-\alpha(m))x_1}\,dx_1 \nonumber\\ 
& & +\ \frac{1}{2R}\int\limits_{-R}^R u^{sc}_{n,pert}(x)\,e^{-i\alpha(m)x_1}\,dx_1\,. \label{eq1}
\end{eqnarray}
We first estimate the first and second terms on the right hand side of the above relation. 
It is obvious that there exists $\delta>0$ such that
$$ \min\bigl\{|\ell+\alpha(n)-\alpha(m)|:n,m\in\{1,\ldots,N\},\ \ell\in\ganz,\ 
\ell+\alpha(n)-\alpha(m)\not=0\bigr\}\ >\ \delta\,. $$ 
In the particular case that $\ell=0$, we have $$|\alpha(n)-\alpha(m)|>\delta\quad\mbox{ for all} 
\quad n\neq m,\ n,m=1,2,\cdots,N\,, $$
because $\alpha(n)\not=\alpha(m)$ for $n\not=m$. We compute explicitly the first and second  integrals as follows:
\begin{eqnarray*}
& & \sum_{n=1}^N\lambda_n\,e^{-i\sqrt{k^2-\alpha(n)^2}x_2}\frac{1}{2R}\int\limits_{-R}^R
e^{i(\alpha(n)-\alpha(m))x_1}\,dx_1 \\
& = & \lambda_m\,e^{-i\sqrt{k^2-\alpha(m)^2}x_2}\ +\sum_{n\neq m, n=1}^N
\lambda_n\,e^{-i\sqrt{k^2-\alpha(n)^2}x_2}\frac{1}{R} \frac{\sin [ \alpha(n)-\alpha(m) R]}
{\alpha(n)-\alpha(m)} \\
& = & \lambda_m\,e^{-i\sqrt{k^2-\alpha(m)^2}x_2}\ +\ \cO(1/R)\,,
\end{eqnarray*}
and
\begin{eqnarray*}
& & \sum_{n=1}^N\sum_{\ell\in\ganz}\lambda_n\,u_{\ell,n}\,e^{i\sqrt{k^2-(\ell+\alpha(n))^2}x_2}\,
\frac{1}{2R}\int\limits_{-R}^R e^{i(\ell+\alpha(n)-\alpha(m))x_1}\,dx_1 \\
& = & \sum_{n=1}^N\sum_{\ell\in\ganz}\lambda_n\,u_{\ell,n}\,e^{i\sqrt{k^2-(\ell+\alpha(n))^2}x_2}\,
\frac{\sin[(\ell+\alpha(n)-\alpha(m))R]}{(\ell+\alpha(n)-\alpha(m))R} \\
& = & e^{i\sqrt{k^2-\alpha(m)^2}x_2}\sum_{n,\ell:\ell+\alpha(n)=\alpha(m)}\lambda_n\,u_{\ell,n}\,\
+\ \cO(1/R).
\end{eqnarray*}
 Next we estimate the third term of (\ref{eq1}). By definition of the open waveguide 
radiation condition, we can decompose $u^{sc}_{n,pert}$ into the sum $u^{sc}_{n,pert}= 
u^{sc,prop}_{n,pert}+u^{sc,rad}_{n,pert}$ in $\tilde{D}$ where the radiating part 
$u^{sc,rad}_{n,pert}\in H^1(W_h)$, and the propagating part $u^{sc,prop}_{n,pert}$ is of the form 
\eqref{u2}. The term involving $u^{sc,rad}_{n,pert}$ converges to zero as $R\rightarrow\infty$, 
because 
\begin{eqnarray*}
\frac{1}{2R}\biggl|\int\limits_{-R}^Ru^{sc,rad}_{n,pert}(x_1,x_2)\,e^{-i\alpha(m)x_1}\,dx_1\biggr| 
& \leq & \frac{1}{2R}\int\limits_{-R}^R\big|u^{sc,rad}_{n,pert}(x_1,x_2)\big|\,dx_1 \\
& \leq & \frac{1}{\sqrt{2R}}\,\Vert u^{sc, rad}_{n,pert}(\cdot,x_2)\Vert_{L^2(\R)}
\end{eqnarray*}
for all $x_2\geq h$. To estimate the propagating part, we observe that for $x_2>h$ it has the 
form
\begin{eqnarray*}
u^{sc,prop}_{n,pert}(x) & = & \psi^+(x_1)\sum_{j\in J}\sum_{\ell:|\ell+\hat{\alpha}_j|>k}
c^+_{\ell,j,n}\,e^{-\sqrt{(\ell+\hat{\alpha}_j)^2-k^2}x_2}\,e^{i(\ell+\hat{\alpha}_j)x_1} \\
& & +\ \psi^-(x_1)\sum_{j\in J}\sum_{\ell:|\ell+\hat{\alpha}_j|>k}c^-_{\ell,j,n}\,
e^{-\sqrt{(\ell+\hat{\alpha}_j)^2-k^2}x_2}\,e^{i(\ell+\hat{\alpha}_j)x_1}
\end{eqnarray*}
for some coefficients $c^\pm_{\ell,j,n}$. Therefore,
\begin{eqnarray*}
& & \frac{1}{2R}\int\limits_{-R}^Ru^{sc,prop}_{n,pert}(x_1,x_2)\,e^{-i\alpha(m)x_1}\,dx_1 \\
& = & \sum_{j\in J}\sum_{\ell:|\ell+\hat{\alpha}_j|>k}c^+_{\ell,j,n}\,
e^{-\sqrt{(\ell+\hat{\alpha}_j)^2-k^2}x_2}\frac{1}{2R}\int\limits_{-R}^R\psi^+(x_1)\,
e^{i(\ell+\hat{\alpha}_j-\alpha(m))x_1}\,dx_1 \\
& & +\ \sum_{j\in J}\sum_{\ell:|\ell+\hat{\alpha}_j|>k}c^-_{\ell,j,n}\,
e^{-\sqrt{(\ell+\hat{\alpha}_j)^2-k^2}x_2}\frac{1}{2R}\int\limits_{-R}^R\psi^-(x_1)\,
e^{i(\ell+\hat{\alpha}_j-\alpha(m))x_1}\,dx_1\,,
\end{eqnarray*}
and 
\begin{eqnarray*}
& & \frac{1}{2R}\int\limits_{-R}^R\psi^+(x_1)\,e^{i(\ell+\hat{\alpha}_j-\alpha(m))x_1}\,dx_1\\
& = & \frac{1}{2R}\int\limits_{\sigma-1}^{\sigma}\psi^+(x_1)\,
e^{i(\ell+\hat{\alpha}_j-\alpha(m))x_1}\,dx_1\ +\
\frac{1}{2R}\int\limits_{\sigma}^Re^{i(\ell+\hat{\alpha}_j-\alpha(m))x_1}\,dx_1
\end{eqnarray*}
converges to zero as $R$ tends to infinity uniformly with respect to $\ell$ by the same arguments 
as in part (i) since $|\ell+\hat{\alpha}_j|>k>|\alpha(m)|$. The same argument applies to the term 
involving $\psi^-(x_1)$.
Letting $R$ tend to infinity in \eqref{eq1}  we conclude that
$$ \lambda_m\,e^{-i\sqrt{k^2-\alpha(m)^2}x_2}\ +\ e^{i\sqrt{k^2-\alpha(m)^2}x_2}\biggl[
\sum_{n,\ell:\ell+\alpha(n)=\alpha(m)}\lambda_n\,u_{\ell,n}\biggr]\ =\ 0\,. $$
The linear independence of the exponential terms yields that $\lambda_m=0$. This ends the proof of the linear independence of the total fields with different directions.

\smalf
Finally, repeating the lines in the proof of Theorem \ref{uni-sour-inf} with $G_1(x;y_j)=u_1(x;\theta_m)$, we can prove the uniqueness by the same contradiction argument. 
\end{proof}
 In the appendix we shall provide another proof of the linear independence of the total fields $\{u(x;\theta_n)\}_{n=1}^N$ for any $N\in \N$.

\subsection{Uniqueness With A Finite Number of Plane Waves}

In Theorems \ref{uni-sour-inf} and \ref{uni-plane-inf}, there is no requirement on the location, width and height of the defect. If some a priori information on the defect is available, we can prove uniqueness with a finite number of incoming waves by adopting Colton and Slemann's idea of determining a bounded sound-soft obstacle \cite{CS83}. 

\begin{theorem}\label{uni-plane-fini} Let $k>0$ be fixed and let $\tilde{\Gamma}$ be a local perturbation of the periodic curve $\Gamma$. Suppose that 
$\max\{x_2: x\in \Gamma\cup \tilde{\Gamma}\}<h$ for some $h\in \R$ and that both $\Gamma\backslash\tilde{\Gamma}$ and $\tilde{\Gamma}\backslash\Gamma$ are contained in the rectangular domain $D_0=(0, 2\pi)\times(0, h)$. Let $N\geq hk^2/2$ be an integer.  
Then $\tilde{\Gamma}$ can be uniquely determined by the near-field measurement data $\{u(x_1, h; \theta_n): x_1\in(a,b),  n=1,2,\cdots, N+1\}$ where  $\theta_n\in (-\pi/2, \pi/2)$ are distinct angles.
\end{theorem}
\begin{proof} Suppose that there are two local perturbations $\tilde{\Gamma}_1$ and $\tilde{\Gamma}_2$ lying below the line $x_2=h$, which produce identical near-field data for each incident direction $\theta_m$. Denote by $u_{\ell}(x;\theta_m)$ ($\ell=1,2$) the unique total field incited by the incoming plane wave $u^{in}(x;\theta_m)=e^{ik (x_1\sin\theta_m -x_2\cos\theta_m  )}$ incident onto  $\tilde{\Gamma}_{\ell}$. 
We proceed as in the proof of Theorem \ref{uni-sour-inf}  to obtain
\ben
(\Delta +k^2)u_1(x; \theta_n)=0\quad\mbox{in}\quad D^*,\quad u_1(\cdot; \theta_n)=0\quad\mbox{on}\quad \partial D^*
\enn for all $n=1,2,\cdots, N, N+1$, where $D^*\subset D_0$ is a bounded domain. This implies that there exist $N+1$ Dirichlet eigenfunctions $u_1(\cdot; \theta_n)\in H_0^1(D^*)$ for the negative Laplacian operator over the bounded domain $D^*$ with the eigenvalue $k^2$. Recalling the linear independence of $u_1(\cdot; \theta_n)$ (see the proof of Theorem \ref{uni-plane-inf}), we conclude that the dimension of the Dirichlet eigenspace over $D^*$ associated with $k^2$ must be greater than or equal to $N+1$. Below we shall prove that this dimension cannot exceed $N$, which is a contradiction.

\smalf
Denote by $\lambda_j$ $(j\in \N)$ the Dirichlet eigenvalues of $D^*$, which are arranged according to increasing magnitude and taken with respective to multiplicity. Let the multiplicity of $k^2$ be  $m^*\in \N$ and suppose that $k^2$ is the $m$-th ($m\geq m^*$) eigenvalue such that 
\ben
\lambda_{m+1}>k^2=\lambda_m=\lambda_{m-1}=\cdots=\lambda_{m-m^*-1}>\lambda_{m-m^*-2}\geq \lambda_{m-m^*-3}\geq \cdots >\lambda_1>0.
\enn
Analogously, let $0<\mu_1\leq \mu_2\leq\cdots\leq \mu_m$ be the first $m$ eigenvalues of $D_0$.
By the strong monotonicity property of the Dirichlet eigenvalues with respect to the domain, it holds that $\mu_m<\lambda_m=k^2$ due to the fact that $D^*\subset D_0$. This further implies that $m^*$ is less than or equal to $(D_0, k^2)^\sharp\in \N$, which is defined as the sum of the multiplicities of the Dirichlet eigenvalues for $D_0$ that are less than $k^2$. On the other, if $k^2$ is a Dirichlet eigenvalue of the rectangular domain $D_0$, it is easy to derive using the method of separating variables the associated eigenfunctions 
\ben
v_{l,j}(x_1, x_2)=\sin\left(\frac{l}{2}x_1\right)\,\sin\left(\frac{j\pi}{h}x_2\right),
\enn
where $l, j\in \N$ satisfy the relation
\be\label{eigenvalue}
k^2=\frac{l^2}{4}+\frac{j^2\pi^2}{h^2}.
\en
Therefore,  $(D_0, k^2)^\sharp$ coincides with the number of grid points $(l,j)\in \N\times\N$ lying in the positive orthant of the ellipse 
\ben
\left(\frac{x_1}{2k}\right)^2+\left(\frac{x_2}{hk/\pi}\right)^2\leq 1.
\enn
Hence, $(D_0, k^2)^\sharp$ can be bounded by $hk^2/2$, one fourth of the area of the above ellipse.   By the choice of $N$, we have $(D_0, k^2)^\sharp\leq N$. This contradicts the fact that there are $N+1$ linearly independent functions $u_1(x;\theta_m)\in H_0^1(D^*)$ for $m=1,2,\cdots, N+1$.
\end{proof}
As a direct consequence of the proof of Theorem \ref{uni-plane-fini}, we can obtain a uniqueness result with one plane wave with fixed direction and frequency.
\begin{corollary}\label{uni-plane-one} Let $k>0$ be fixed and let $\tilde{\Gamma}$ be a local perturbation of the periodic curve $\Gamma$. Suppose that 
$\max\{x_2: x\in \Gamma\cup \tilde{\Gamma}\}<h$ for some $h\in \R$ and that both $\Gamma\backslash\tilde{\Gamma}$ and $\tilde{\Gamma}\backslash\Gamma$ are contained in the rectangular domain $D_0=(0, 2\pi)\times(0, h)$. If $k<\sqrt{1/4+\pi^2/h^2}$,
then $\tilde{\Gamma}$ can be uniquely determined by a single near-field measurement data $\{u(x_1, h; \theta): x_1\in(a,b)\}$ where  $\theta\in (-\pi/2, \pi/2)$ is arbitrary.
\end{corollary}
\begin{proof} From the proof of Theorem \ref{uni-plane-fini}, we conclude that $k^2$ must be greater than or equal to the first Dirichlet eigenvalue $\mu_1$ of the negative Laplacian operator over the domain $D_0$. In view of \eqref{eigenvalue}, one obtains $\mu_1= \sqrt{1/4+\pi^2/h^2}\leq k^2$, which is a contradiction to the condition that $k<\sqrt{1/4+\pi^2/h^2}$.
Hence, $u(x;\theta)$ cannot be a Dirichlet eigenfunction over any subdomain of $D_0$. This proves the desired uniqueness result by applying the same contradiction arguments of Theorem \ref{uni-plane-fini}. 
\end{proof}

\smalf

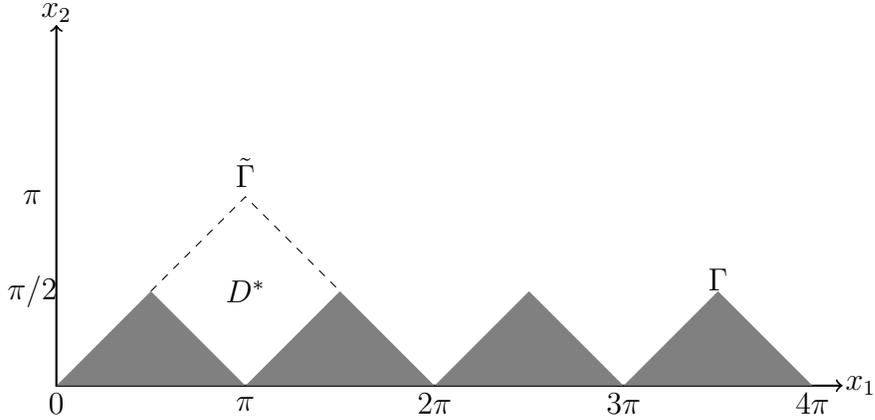
\begin{figure}[htb]
\centering

\begin{tikzpicture}[scale=.8, thick]

\draw [<->] (0,6) -- (0,0) -- (4*pi+0.5,0);
\fill[gray] (0,0) -- (pi/2, pi/2 ) -- (pi, 0) -- (3*pi/2, pi/2) -- (2*pi, 0)
                     -- (5*pi/2, pi/2 ) -- (3*pi, 0) -- (7*pi/2, pi/2) -- (4*pi, 0);

\draw [thin, dashed]  (pi/2, pi/2 ) --(pi, pi)--(3*pi/2, pi/2);

\draw (0,-0.3) node{$0$};
\draw (pi,-0.3) node{$\pi$};
\draw (2*pi,-0.3) node{$2\pi$};
\draw (3*pi,-0.3) node{$3\pi$};
\draw (4*pi,-0.3) node{$4\pi$};
\draw (-0.4,pi/2) node{$\pi/2$};
\draw (-0.4,pi) node{$\pi$};
\draw (4*pi+0.8,0) node{$x_1$};
\draw (0, 6.2) node{$x_2$};

\draw (pi, pi/2) node{$D^*$};
\draw (7*pi/2, pi/2+0.2) node{$\Gamma$};
\draw (pi, pi+0.4) node{$\tilde{\Gamma}$};
\end{tikzpicture}
\caption{Illustration of $\Gamma$ and its local perturbation $\tilde{\Gamma}$ which generate identical wave fields with $\theta=0$ and $k=2$. Here $D^*=D\backslash\overline{\tilde{D}}$ represents the difference domain between $D$ and $\tilde{D}$.}\label{f2}
\end{figure}

Below we present a counterexample to show that, if $k\geq\sqrt{1/4+\pi^2/h^2}$, it is in general impossible to unique determine the defect using a single plane wave when the Rayleigh frequency occurs. Such an example is motivated by the classification of unidentifiable polygonal gratings with one acoustic or elastic plane wave \cite{EH10, EH11}.
Let the incident angle be $\theta=0$ and set $k=2$, leading to $u^{in}(x;\theta)=e^{-i2x_2}$. Define the piecewise linear function (see Figure \ref{f2})
\ben
x_2=f(x_1)=\left\{\begin{array}{lll}
x_1&&\quad\mbox{if}\quad x_1\in(0, \pi/2),\\
-x_1+\pi &&\quad\mbox{if}\quad x_1\in(\pi/2, \pi).
\end{array}\right.
\enn
Let $\Gamma$ be the $\pi$-periodic extensions of $\{x_2=f(x_1): x_1\in(0, \pi)\}$ in the $x_1$-direction, and let $\tilde{\Gamma}$ be the local perturbation of $\Gamma$ in $(0, 2\pi)$ shown as in Figure \ref{f2}, where the dashed line segments denote the defect and $D^*$ the gap domain between $D$ and $\tilde{D}$. In this case we have $k>\sqrt{1/4+\pi^2/h^2}$ for all $h>\pi$. Since $\Gamma$ is the graph of a piecewise linear function, there exists a unique scattered field $u_{unpert}^{sc}\in H^1_{loc,\alpha}(D)$  to the unperturbed scattering problem,  taking the explicit form
\ben
u_{unpert}^{sc}(x)=e^{i2x_2}-e^{i2x_1}-e^{-i2x_1},\quad x\in D.
\enn
Note that the Rayleigh frequency occurs, since $k=2$ and $\alpha=0$ (that is, $k=\alpha+n$ with $n=2$).  Moreover, the guide modes (surface waves) are excluded for the unperturbed scattering problem. Hence, the unique total field to the unperturbed problem can be expressed as
\ben
u_{unpert}=u^{in}+u_{unpert}^{sc}=2\left(\cos 2x_2-\cos 2 x_1\right)=4\sin(x_1+x_2)\,\sin(x_1-x_2),\; x\in D.
\enn
Observing that $\tilde{\Gamma}$ is also the graph of a piecewise linear function, 
by \cite{HWR} there exists a unique total field $u_{pert}\in H^1_{loc, 0}(\tilde{D})$  of the form $u_{pert}=u_{unpert}+u^{sc}_{pert}$ in $\Sigma_R$, where $u^{sc}_{pert}\in H^1_{loc}(\tilde{D})$ consists of the radiating part only satisfying the Sommerfeld radiation conditions of Definitions \ref{src} and \ref{src-i}. On the other hand, since the defect lies on the straight lines 
\ben
x_2=x_1,\quad x_2=-x_1+2\pi, \quad x_1\in \R,
\enn
we conclude that $u_{unperp}$ also vanishes on $\tilde{\Gamma}$. By uniqueness, this implies that $u^{sc}_{perp}\equiv 0$ in $\tilde{D}$ and thus $u_{perp}= u_{unperp}$. In other words, the presence of the local defect does not produce any perturbation to $u_{unperp}$. We remark that $u_{unpert}$ is a real-valued Dirichlet eigenfunction of the negative Laplacian operator over the rectangular domain $D^*$.
Therefore, it is impossible to determine the perturbed boundary $\tilde{\Gamma}\backslash\Gamma$ from the near-field measurement data of $u_{perp}$ on $x_2=h$ for all $h>\pi$. 

\subsection{Uniqueness Using Far-field Data of Point Source Waves}\label{sub5-4}
The symmetry of the Green's function (see Theorem \ref{sy}) together with Theorem \ref{TH-PL} yields
\be\nonumber
u^\infty (-\hat{\theta}; x)&:=& \lim_{t\rightarrow\infty} \left[\sqrt{t} e^{-ikt} u(z_t;x)\right] \\ \nonumber
&=&\lim_{t\rightarrow\infty} \left[\sqrt{t} e^{-ikt} u(x; z_t)\right]\\ \label{mix}
&=&\gamma\,w(x; \hat{\theta})
\en
in $H^1(\tilde{D}_R)$  for any $R>\pi$ with $\gamma:=\frac{e^{i\pi/4}}{\sqrt{8k\pi}}$. Here $z_t=t\hat{\theta}$ and $w(x;\hat{\theta})$ denotes the total field excited by the plane wave $u^{in}(x;\theta)=e^{ikx\cdot\hat{\theta}}$ (see Theorem \ref{TH-PL}). This means that
 the far-field data at the observation direction $-\hat{\theta}$ generated by the point source wave emitting from $x\in \tilde{D}$ is identical with the value of the total field at $x$ of the plane wave $e^{ik x\cdot \hat{\theta}}$ multiplied by the constant $\gamma$. This is exactly the mixed reciprocity relation between point source and plane wave incidences in a perturbed periodic structure. 
Using this relation we can prove uniqueness in determining the defect from the far-field patterns of one or many point source waves. 
\begin{theorem}\label{uni-point-inf} Let $\tilde{\Gamma}$ be a local perturbation of the periodic curve $\Gamma$ and suppose that 
$\max\{x_2: x\in \Gamma\cup \tilde{\Gamma}\}<h$ for some $h\in \R$.  
Then $\tilde{\Gamma}$ can be uniquely determined by the far-field measurement data $\{u^\infty(-\hat{\theta}_m; x): m=1,2,\cdots, x=(x_1, h), x_1\in(a, b)\}$ incited by infinitely many point waves lying on the line segment $\{(x_1, h): x_1\in(a, b)\}$. The same uniqueness result holds true if we replace $u^\infty(-\hat{\theta}; x)$ by $u_{perp}^\infty(-\hat{\theta}; x)$, the far-field pattern of the radiating part of the scattered field $u^{sc}_{perp}(t\hat{\theta}; x)$ as $t\rightarrow\infty$.
\end{theorem}
\begin{proof}
The first assertion  follows directly from the the mixed reciprocity relation \eqref{mix} and the uniqueness result of Theorem \ref{uni-plane-inf}. To prove the second assertion, we recall from Proposition \ref{wps} a decomposition of $u(x; z_t)$ into 
\be\label{dec1}
u(x; z_t)=u^{in}(x; z_t)+u^{sc}_{unperp}(x; z_t)+u^{sc}_{perp}(x; z_t)\quad\mbox{in}\quad \Sigma_R,
\en
where  $u^{sc}_{unperp}(x; z_t)$ denotes the scattered field to the unperturbed scattering problem and $u^{sc}_{perp}(x; z_t)$ the part cased by the defect. Note that both $u^{sc}_{unperp}(x; z_t)$ and $u^{sc}_{perp}(x; z_t)$ fulfil the open waveguide radiation condition. Since the propagating part of $u^{sc}_{unperp}(x; z_t)$ (resp. $u^{sc}_{perp}(x; z_t)$) decays exponential as $t\rightarrow\infty$, one deduces from \eqref{dec1} the corresponding  decomposition of the far-field pattern:
\ben
u^\infty(\hat{\theta}; x)=e^{ik\hat{\theta}\cdot x}+u_{unperp}^\infty(\hat{\theta}; x)+u_{perp}^\infty(\hat{\theta}; x),
\enn
where $u_{unperp}^\infty(\hat{\theta}; x)$ represents the far-field pattern of the radiating part of the scattered field $u^{sc}_{unperp}(t\hat{\theta}; x)$ as $t\rightarrow\infty$. 
Since the unperturbed structure $\Gamma$ is a priori given, 
the far-field pattern $u_{unperp}^\infty(\hat{\theta}; x)$ is uniquely determined by $e^{ik\hat{\theta}\cdot x}$ and $\Gamma$. Hence, the knowledge of  $u_{perp}^\infty(\hat{\theta}; x)$ is equivalent to knowing $u^\infty(\hat{\theta}; x)$ for any fixed $\hat{\theta}\in \mathbb{S}$ and  $x\in \tilde{D}$. This proves the second assertion of Theorem \ref{uni-point-inf}.
\end{proof}

If {\rm a priori} information on the height and size of the defect is available, one can also determine the defect by taking the far-field measurement data at a finite number of observation directions excited by infinitely many point source waves. 

\begin{corollary} Let the conditions of Theorem \ref{uni-plane-fini} hold. Then $\tilde{\Gamma}$ can be uniquely determined by the far-field data $u^\infty(\hat{\theta}_n; x_1, h)$ (or $u^\infty_{perp}(\hat{\theta}_n; x_1, h) )$  for all $x_1\in(a,b)$,  $n=1,2,\cdots, N+1$ where  $\hat{\theta}_n\in \mathbb{S}_+:=\{(x_1, x_2)\in \mathbb{S}: x_2>0\}$ are distinct observation directions. Moreover, the far-field data at a single observation direction (i.e., $N=0$) are sufficient under the additional conditions of Corollary \ref{uni-plane-one}.   
\end{corollary}

\begin{remark}
It is unclear to us the uniqueness with far-field patterns of plane wave incidences, due to the lack of the one-to-one correspondence between the far-field pattern and near-field data. 
\end{remark}

\section{Appendix}
\subsection{An Alternative Proof To The Linear Independence Of Total Fields With Different Directions}
Here we present another proof by adopting the arguments of \cite[Lemma 2.2]{XHZZ}.
Let $u(x; \theta_m)$ with $\theta_m\in(-\pi/2, \pi/2)$ be the uniquely determined total field of the perturbed scattering problem; see Theorems \ref{TH-LAP} and \ref{TH-POI}.  Suppose that $\sum_{m=1}^M c_m u(x; \theta_m)=0$ for all $x\in \tilde{D}$, where $c_m\in \C$, $m=1,2,\cdots, M$ with some $M\in \N$. By \eqref{total},
\ben
\sum_{m=1}^M c_m u_{unpert}(x; \theta_m)+\sum_{m=1}^M c_m u^{sc}_{pert}(x; \theta_m)=0\quad\mbox{for all}\quad x\in U_h,\quad h>h_0.
\enn
In view of the definition of the open waveguide radiation condition, $u^{sc}_{unpert}$ can be decomposed into two parts,
\be\label{dec}
u^{sc}_{unpert}(x; \theta_m)=u^{rad}_{unpert}(x; \theta_m)+u^{prop}_{unpert}(x; \theta_m)\quad \mbox{in}\quad \Sigma_R 
\en
 where  the radiating part $u^{rad}_{unpert}(x; \theta_m)$ decays as $|x|^{-1/2}$ in $U_h$ as $|x|\rightarrow\infty$, whereas the propagating part $u^{prop}_{unpert}(x; \theta_m)$ exponentially decays as $x_2\rightarrow\infty$. This leads to 
\ben
\lim_{H\rightarrow\infty}\frac{1}{H}\Big|\int\limits_{H}^{2H} u^{sc}_{pert}(x; \theta_m)e^{ik\cos\theta_n x_2}\, dx_2\Big|\leq  
\lim_{H\rightarrow\infty}\frac{1}{H}\int\limits_{H}^{2H} \Big|u^{sc}_{pert}(x; \theta_m)\Big|\, dx_2=0
\enn
for all $n, m=1,2,\cdots, M$ and uniformly in all $x\in[0, 2\pi]$. Now, multiplying $e^{ik\cos\theta_n x_2}$ with some $n\in\{1,2,\cdots, M\}$ to both sides of \eqref{dec}, integrating over $(H, 2H)$ with respect to $x_2$ and taking the limit as $H\rightarrow \infty$ yields
\ben
\lim_{H\rightarrow\infty}\frac{1}{H}\Big|\int\limits_{H}^{2H} \sum_{m=1}^M c_m\, u_{unpert}(x; \theta_m)e^{ik\cos\theta_n x_2}\, dx_2\Big|=0\quad\mbox{for all}\quad x_1\in(0, 2\pi).
\enn
Since the unperturbed total field $u_{unpert}(x; \theta_m)$ is $k\sin\theta_m$-quasiperiodic,   the above relation also holds for all $x_1\in \R$. Therefore, 
\ben
\lim_{H\rightarrow\infty}\frac{1}{H}\Big|\int\limits_{H}^{2H} \sum_{m=1}^M c_m\, \left(e^{ik\sin\theta_m x_1}\, + \, u^{sc}_{unpert}(x; \theta_m)\right) e^{ik\cos\theta_n x_2}\, dx_2\Big|=0,\qquad x_1\in \R.
\enn
By the proof of \cite[Lemma 2.2]{XHZZ}, the previous relation implies $c_n=0$. Using the arbitrariness of $1\leq n\leq M$, one obtains $c_m=0$ for all $m=1,2,\cdots, M$. This proves that $\{u(x;\theta_m)\}_{m=1}^M$ must be linearly independent for all $M\in \N$.

\subsection{Proof Of The Asymptotics In \eqref{partial}. }\label{asy}
We suppose that $I_{\tilde{m}}=(a, b)\subset [-1/2, 1/2]$, where the two ending points $a$ and $b$ maybe the cut-off values of $k$. For any $l\in \mathcal{L}\backslash\{\tilde{l}\}$, it holds that $k>|l+\alpha|$. Set $f(\alpha):=(l+\alpha)\sin\theta+\sqrt{k^2-(l+\alpha)^2}\cos\theta$ for $\alpha\in(a, b)$. It is easy to see 
\ben
f'(\alpha)=\sin\theta-\frac{l+\alpha}{\sqrt{k^2-(l+\alpha)^2}}\cos\theta\neq 0\quad\mbox{for all}\quad
l\in \mathcal{L}\backslash\{\tilde{l}\}, \alpha\in (a, b),
\enn and that
\ben
g(\alpha):=\sqrt{k^2-(l+\alpha)^2}\, f'(\alpha)=\sqrt{k^2-(l+\alpha)^2}\sin\theta-(l+\alpha)\cos\theta
\enn
keeps a positive distance from zero for all $
l\in \mathcal{L}\backslash\{\tilde{l}\}$ and $\alpha\in [a, b]$. Here and below the prime always denotes the derivative with respect to $\alpha$.
 Direct calculations show that
\be\nonumber
&&\int\limits_a^b \frac{e^{i t f(\alpha)}\,v_{l, \alpha}}{\sqrt{k^2-(l+\alpha)^2}} \,d \alpha=\frac{-i}{t}\int\limits_a^b
\frac{v_{l, \alpha}}{g(\alpha)} \;d\, e^{i t f(\alpha)}\\ \label{a1}
&=&\frac{-i}{t} \left[   \frac{e^{i t f(\alpha)}\;v_{l, \alpha}}{g(\alpha)}\Big|_{a}^b- \int\limits_a^b
\left( \frac{v_{l, \alpha}'}{g(\alpha)} -v_{l, \alpha} \frac{g'(\alpha)}{g^2(\alpha)} \right) e^{i t f(\alpha)} \,d\alpha                  \right].
\en
We note that $g(\alpha)$ does not vanish at the boundary points $a$ and $b$ for $l\neq \tilde{l}$ and that 
\ben
\int\limits_a^b |g'(\alpha)| d\alpha=\int\limits_a^b\Big |\frac{l+\alpha}{\sqrt{k^2-(l+\alpha)^2}}\sin\theta+\cos\theta\Big| d\alpha<\infty,
\enn
because $g'$ has integrable singularities at the possible cut-off values on $a$ or $b$. 
On the other hand,
$v_{l, \alpha}\in H^1_{\alpha, loc, 0}(D)$ can be chosen to depend continuously on $\alpha\in[-1/2, 1/2]$ (see \cite[Theorem 3.3]{HK22}) and $v_{l, \alpha}\in W^{1,1}([-1/2, 1/2]; H^1(Q_{h_0}))$ has only square-root singularities at cut-off values, which can be verified following the same arguments in the proof of \cite[Theorem 4.3]{K22} for inhomogeneous layers. Hence, the right hand side of \eqref{a1} decays as $O(t^{-1})$ as $t$ tends to infinity. This together with the arbitrariness of $l\in \mathcal{L}\backslash\{\tilde{l}\}$ proves the first relation in \eqref{partial}. The second one can be verified analogously by noting that $f'(\alpha)\neq 0$ for all $l\in \mathcal{L}$ and $\alpha\in [-1/2, 1/2]\backslash (a, b)$.


\section*{Acknowledgements}
The first author (G.H.) acknowledges the hospitality of the Institute for Applied and Numerical 
Mathematics, Karlsruhe Institute of Technology and the support of Alexander von 
Humboldt-Stiftung.  The second author (A.K.) gratefully acknowledges the financial 
support by the Deutsche Forschungsgemeinschaft (DFG, German Research Foundation) -- 
Project-ID 258734477 -- SFB 1173.


\begin{thebibliography}{10}

\bibitem{Lu19} A. Abdrabou and Y. Lu, Indirect link between resonant and guided modes on uniform and periodic slabs, Physical Review A. 99 (2019): 063818.

\bibitem{Bao2022}
G.~Bao and P.~Li.
\newblock {\em Maxwell's equations in periodic structures}, volume 208 of {\em
  Applied Mathematical Sciences}.
\newblock Springer, Singapore; Science Press Beijing, Beijing, 2022.

\bibitem{BBS94} A. S. Bonnet-Bendhia and P. Starling, Guided waves by electromagnetic 
gratings and non-uniqueness examples for the diffraction problem, Math. Meth. Appl. Sci., 
17 (1994):  2305-338.


\bibitem{CE10} S. N. Chandlea-Wilde and J. Elschner,  Variational approach in weighted 
Sobolev spaces to scattering by unbounded rough surfaces, SIAM J. Math. Anal., 42 (2010): 
2554-2580.

\bibitem{SM05}S. N. Chandlea-Wilde and P. Monk,  Existence, uniqueness, and variational 
methods for scattering by unbounded rough surfaces, SIAM J Math. Anal., 37 (2005): 598-618.

\bibitem{CR1995} S. N. Chandler-Wilde and C.R. Ross,  Uniqueness results for direct and inverse scattering by infinite surfaces in a lossy medium, Inverse Problems 11 (1995): 1063-1067.


\bibitem{CZ98} S.N. Chandler-Wilde and B. Zhang,  A uniqueness result for scattering by infinite rough surfaces,
SIAM J. Appl. Math., 58 (1998): 1774-1790.



\bibitem{CK13} D. Colton and R. Kress, {\rm Integral Equation Methods in Scattering Theory}, Volume 72 of Classics in Applied Mathematics, Society for Industrial and Applied Mathematics, 2013. 

\bibitem{CS83} D. Colton and B. D. Sleeman,   Uniqueness theorems for the inverse problem of acoustic scattering,  IMA J. Appl. Math. 31 (1983): 253-259. 


\bibitem{EH10} J. Elschner and G. Hu, Global uniqueness in determining polygonal periodic structures with a minimal number of incident plane waves, Inverse Problems, 26 (2010):  115002/1--115002/23.

\bibitem{EH11} J. Elschner and G. Hu, Inverse scattering of elastic waves by periodic structures: Uniqueness under the third or fourth kind boundary conditions, Methods and Applications of Analysis 18 (2011): 215-244.


\bibitem{EY02} J. Elschner and M. Yamamoto, An inverse problem in periodic diffractive optics: 
Reconstruction of Lipschitz grating profiles, Appl. Anal., 81 (2002):  1307-1328.


\bibitem{FJ09} S. Fliss and P. Joly, Exact boundary conditions for time-harmonic wave 
propagation in locally perturbed periodic media, Appl. Numer. Math., 59 (2009): 2155-2178.

\bibitem{FJ16} S. Fliss and P. Joly, Solutions of the time-harmonic wave equation in periodic 
waveguides: Asymptotic behavior and radiation condition, Arch. Ration. Mech. Anal., 219 
(2016):  349-386.

\bibitem{TF} T. Furuya, Scattering by the local perturbation of an open periodic waveguide 
in the half plane, J. Math. Anal. Appl., 489 (2020): 124-149.

\bibitem{G00} V. Yu. Gotlib, Solutions of the Helmholtz equation, concentrated near a plane 
periodic boundary, J. Math. Sci., 102 (2000): 4188-4194.

\bibitem{Gri16} A. Grigor'yan,   Analysis of Elliptic Differential Equations, A Lectrure note in 2016, Universität Bielefeld, Germany.  https://www.math.uni-bielefeld.de/~grigor/elelect.pdf

\bibitem{HF97} F. Hettlich and A. Kirsch, Schiffer’s theorem in inverse scattering theory for periodic structures,
Inverse Problems, 13 (1997): 351–361. 


\bibitem{HK22} G. Hu and A. Kirsch, Time-harmonic scattering by locally perturbed periodic structures with Dirichlet and Neumann boundary conditions, arXiv 2024.

\bibitem{HWR} G. Hu, W. Lu and A.Rathsfeld, Time-harmonic acoustic scattering from locally 
perturbed periodic curves, SIAM J. Appl. Math., 81 (2021): 2569-2595.

\bibitem{HR15}
G. Hu and A. Rathsfeld, Scattering of time-harmonic electromagnetic plane waves by perfectly conducting diffraction gratings, IMA Appl. Math. 80 (2015): 508-532.

\bibitem{Isa90} V. Isakov, On uniqueness in the inverse transmission scattering problem, Comm. Part. Diff.
Equat., 15 (1990): 1565--1587.

\bibitem{JLF2006} P. Joly, J.-R. Li and S. Fliss, Exact boundary conditions for periodic 
waveguides containing a local perturbation, Commun. Comput. Phys., 1 (2006): 945–973. 

\bibitem{KN02} I. V.  Kamotski and  S. A. Nazarov, The augmented scattering matrix and 
exponentially decaying solutions of an elliptic problem in a cylindrical domain, J. Math. 
Sci., 111 (2002): 3657-3666.

\bibitem{K93} A. Kirsch, Diffraction by periodic structures, in Proceedings of the Lapland 
Conference on Inverse Problems, L. Paivarinta and E. Summersalo, eds., Springer, Berlin, 
1993, pp. 87- 102.


\bibitem{K19-1} A.Kirsch, Scattering by a periodic tube in  $\mathbb{R}^3$: Part I. The 
limiting absorption principle, Inverse Problems, 35 (2019):  104004 

\bibitem{K19-2} A. Kirsch, Scattering by a periodic tube in  $\mathbb{R}^3$: Part II. A 
radiation condition, Inverse Problems, 35 (2019): 104005.

\bibitem{K22p} A. Kirsch, On the scattering of a plane wave by a perturbed open periodic waveguide, Math. Math. Appl. Sci. 46 (2023): 10698-10718. 

\bibitem{KL18} A. Kirsch and A. Lechleiter, The limiting absorption principle and a radiation 
condition for the scattering by a periodic layer,  SIAM J. Math. Anal.,  50 (2018): 2536-2565.

\bibitem{K22} A. Kirsch, A scattering problem for a local perturbation of an open periodic 
waveguide,  Math. Meth. Appl. Sci., 45 (2022): 5737-5773.

\bibitem{KL-MMAS} A. Kirsch and A. Lechleiter, A radiation condition arising from the 
limiting absorption principle for a closed full- or half-waveguide problem, Math. Meth. Appl. 
Sci.,  41 (2018): 3955--3975.

\bibitem{KK93} A. Kirsch and R. Kress, Uniqueness in inverse obstacle scattering, Inverse Problems, 9 (1993):
285--299.

\bibitem{Kry96} N. V. Krylov, Lectures on Elliptic and Parabolic Equations in H\"older spaces,
American Mathematical Society, 1996. 


\bibitem{LZ17} A. Lechleiter and R. Zhang, A converg nt numerical scheme for scattering of 
aperiodic waves from periodic surfaces based on the Floquet-Bloch transform, SIAM J. Numer. 
Anal., 55 (2017):  713-736.

\bibitem{L2017} A. Lechleiter, The Floquet-Bloch transform and scattering from locally 
perturbed periodic surfaces, J. Math. Anal. Appl., 446 (2017): 605-627.

\bibitem{R07} J. W. S. Lord Rayleigh, On the dynamical theory of gratings, Proc. Roy. Soc. 
Lond. A, 79 (1907): 399-416.



 



\bibitem{M10} W. Mclean, Strongly Elliptic Systems and Boundary Integral Equations, Cambridge 
University Press, Cambridge, UK, 2010.

\bibitem{P1980} R. Petit (ed.). Electromagnetic Theory of Gratings. Springer: Berlin, 1980.




\bibitem{XHZZ} X. Xu, G. Hu, B. Zhang and H. Zhang, Uniqueness to inverse grating diffraction problems with infinitely many plane waves at a fixed frequency, 	SIAM J. Appl. Math 83 (2023):302-326.

\bibitem{Lu18} L. Yuan and Y. Lu, Bound satetes in the continnum on periodic structures surrounded by strong resonances, Physical Review A. 97 (2018): 043828.


\bibitem{XHLR} X. Yu, G. Hu, W. Lu and A. Rathsfeld,
PML and high-accuracy boundary integral equation solver for wave scattering by a locally defected periodic surface,   SIAM Nume.  Anal. 60 (2022): 2592-2625. 

\bibitem{ZR} R. Zhang, Numerical methods for scattering problems in periodic waveguides, Numer. Math. 148 (2021), no. 4, 959–996.



\end{thebibliography}
\end{document}